\documentclass[a4paper,12pt,wlayout]{article}

\usepackage{array} 
\usepackage{amssymb}
\usepackage{amsmath} 
\usepackage{amsbsy} 
\usepackage{epsfig}
\usepackage{mathrsfs} 
\usepackage{graphics}
\usepackage{graphicx}
\usepackage[sort&compress,numbers]{natbib}
\usepackage{pstricks}
\usepackage{float}
\usepackage{amsthm}
\usepackage{slashed}
\usepackage{color}

\theoremstyle{plain}

\theoremstyle{remark}

\def\R{\mathrm{Re}}
\newcommand{\eps}{\varepsilon}
\newcommand{\lb}{\label}
\newcommand{\go}{\rightarrow}

\newcommand{\pt}{\partial_t}
\newcommand{\pxx}{\partial_{xx}}

\newcommand{\pxxx}{\partial_{xxx}}
\newcommand{\px}{\partial_x}
\newcommand{\pzz}{\partial_{zz}}
\newcommand{\pzzz}{\partial_{zzz}}
\newcommand{\pz}{\partial_z}

\newcommand{\ee}{\end{equation}}
\newcommand{\be}{\begin{equation}}
\newcommand{\bea}{\begin{eqnarray}}
\newcommand{\eea}{\end{eqnarray}}
\newcommand{\sbea}{\begin{subequations}\begin{eqnarray}}
\newcommand{\seea}{\end{eqnarray}\end{subequations}} 
\newcommand{\ees}{\end{equation*}}
\newcommand{\bes}{\begin{equation*}}
\newcommand{\beas}{\begin{eqnarray*}}
\newcommand{\eeas}{\end{eqnarray*}}

\newcommand{\rf}[1]{(\ref{#1})}
\newcommand{\myand}{\quad\mbox{and}\quad}
\newcommand{\myat}{\quad\mbox{at}\quad}
\newcommand{\const}{\mathrm{const}}
\newcommand{\mR}{\mathbb{R}\,}
\newcommand{\N}{\mathbb{N}\,}

\usepackage{a4wide}
\begin{document}

\title{Coarsening rates for the dynamics of slipping droplets}
\author{Georgy Kitavtsev\footnote{Max Planck Institute for Mathematics in the Sciences, Inselstr. 22, D--04103 Leipzig, Germany. {\tt E-mail}: Georgy.Kitavtsev@mis.mpg.de}}
\date{\today}
\maketitle
\numberwithin{equation}{section}
\begin{abstract}
We derive reduced finite dimensional ODE models starting from 
one dimensional lubrication equations describing coarsening dynamics of droplets in nanometric polymer film
interacting on a hydrophobically coated solid substrate in the presence of
large slippage at the liquid/solid interface. In the limiting case
of infinite slip length corresponding in applications to free films a collision/absorption model then arises and is solved explicitly. The exact
coarsening law is derived for it analytically and confirmed numerically.
Existence of a threshold for the decay of initial distributions of droplet distances at infinity at which the coarsening rates switch from
algebraic to exponential ones is shown.  
\end{abstract}

\section{Introduction.}
Dewetting processes of a liquid polymer film of nanometer thickness interacting on a hydrophobically coated solid substrate
attracted an intensive research during last several decades, see. e.g. a review in~\cite{ODB97}.
In general, such processes can be divided into three stages.
During the first stage a liquid polymer film is susceptible to instability due to small
perturbations of the film profile. Typically such films rupture,
thereby initiating a complex dewetting process, see e.g.~\cite{Reiter99,RBR91,SHJ01b}. 
The influence of intermolecular forces play an
important part in the rupture and subsequent dewetting process, see
e.g.~\cite{gennes85,WD82} and references therein. 
Typically the competition between the long-range attractive van der Waals and short-range Born
repulsive intermolecular forces reduces the unstable film to \emph{an ultra-thin
layer} that connects the evolving patterns and 
is given by the  minimum of the corresponding intermolecular 
potential, i.e. the film settles into an energetically more favorable state, see 
\cite{EG97,BGW01}.
The second stage is associated with the formation of
regions of this minimal thickness, bounded by moving rims that connect
to the undisturbed film, see e.g.~\cite{SR96,BR92,MW05}. 

\begin{figure}
\centering
\includegraphics*[width=0.4\textwidth]{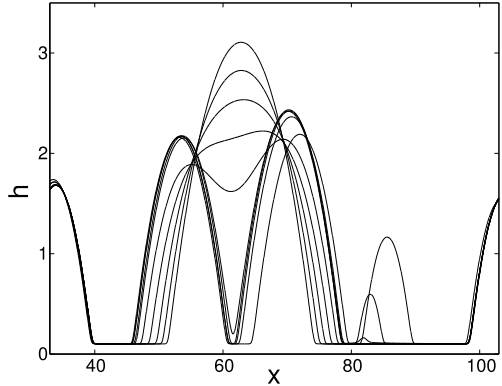}
\caption[Numerical example of coarsening process]
{{\small Numerical solution to \rf{SSM1}--\rf{SSM2} with $\eps=0.1,\,\beta=2.5,\,\R=1$ showing an example
 of a coarsening process (collapse of the 4th small droplet and collision of 2nd and 3rd ones) in the array of five quasiequilibrium droplets.}}
\label{CoarsDyn}
\end{figure}

In this study we are interested in the third and the last stage of the
dewetting process, namely the long-time coarsening process that originates in the breaking up of the evolving patterns into
small droplets and is characterized by its subsequent slow-time coarsening
dynamics, which has been  observed and investigated experimentally by \citet{LG02,LG03}. They show that during the coarsening the average size of droplets increases and the number of droplets decreases. The coarsening mechanisms that were observed in
such films are typically subsequent collapses of smaller droplets and collisions of
neighboring ones. During collapse the size of a droplet shrinks in time and
its mass is distributed in the ultra-thin layer. In turn, collisions among droplets occur due to the mass transfer through the
ultra-thin layer between them that causes a translation movement of them,
\emph{droplet migration}, eventually leading to the formation of new droplets. A numerical example of the coarsening dynamics in two-dimensional films
is shown in Fig. \ref{CoarsDyn}. 
\begin{figure}[ht]
\centering
\includegraphics*[width=0.4\textwidth]{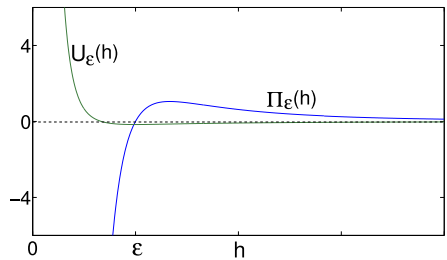}
\caption[Plots of intermolecular pressure and potential function]
{{\small Plots of intermolecular pressure $\Pi_\eps(h)$ (blue) and potential function $U_\eps(h)$ (green) for $\eps=0.1$}}
\label{PFF}
\end{figure}

Besides intermolecular forces and surface tension at the free surface of the film the dewetting of polymer films on hydrophobic substrates also involves such boundary effect as slippage on a solid substrate~\cite{FMWJ07}. Recently in~\citet{MWW06} closed-form one-dimensional lubrication equations over a wide range of slip lengths were derived from the underlying equations for conservation of mass and momentum, together with boundary conditions for the tangential and normal stresses, as well as the kinematic condition at the free boundary, impermeability and Navier-slip condition at the liquid-solid interface. Asymptotic arguments, based on the magnitude of the slip length show that within a lubrication scaling there are two \emph{distinguished regimes}, see~\cite{MWW06}.

These are the well-known \emph{weak-slip} model 
\be
\pt h= - \px \Big(M(h)\px \left(\sigma \pxx h-\Pi_\eps(h)\right)\Big) 
\label{WSM}
\ee
with $M(h):=h^3+b\,h^2$ and $b$ denoting the slip length parameter; and the \emph{strong-slip} model 
\begin{subequations}
\label{SSM}
\begin{align}
\R\,(\eps\pt(hu)+\px(hu^2))&=\nu\px(h\px u)+h\px(\sigma\pxx h-\Pi_\eps(h))-\frac{u}{\beta}
\label{SSM1}\\
\eps\pt h &= -\,\px\left(h u\right),
\label{SSM2}
\end{align}
\end{subequations}
respectively. Here, $u(x,t)$ and $h(x,t)$ denote the average velocity in the
lateral direction and the height profile for the free surface,
respectively. The positive slip length parameters $b$ and $\beta$ are related
by orders of magnitude via $b\sim\eta^2\beta$, where the (small) parameter
$\eta$, $0<\eta\ll 1$, refers to the vertical to horizontal scale separation
of the thin film. 

The high order of the lubrication equations \rf{WSM} and \rf{SSM1}--\rf{SSM2}
is a result of the contribution from surface tension at the free boundary, 
reflected by the linearized curvature term $\sigma \pxx h$ with parameter
$\sigma\ge 0$.  
A further contribution to the pressure is denoted by $\Pi_\eps(h)$ and
represents that of the intermolecular forces, 
namely long-range attractive van der Waals and short-range Born repulsive
intermolecular forces. 
A commonly used expression for it~\cite{BGW01,GW03} is given by
\be
\label{PF} \Pi_\eps(h)=\frac{\eps^2}{h^3}-\frac{\eps^2}{h^4}\ \ \text{with}\ \ 0\le\eps\ll 1. 
\ee
It can be written as a derivative of the potential function
$U_\eps(h)=\mathcal{U}(h/\eps)$ (see Fig. \ref{PFF}) where 
\be
\mathcal{U}(H)=-\frac{1}{2\,H^2}+\frac{1}{3\,H^3},
\label{U(h)}
\ee
The parameter $0<\eps\ll 1$ is the global minimum of $U_\eps(h)$ and
gives to the leading order thickness of the ultra-thin layer. Below we often use notation for the pressure and flux functions function
\be
p(h):=\sigma\pxx h-\Pi_\eps(h),\quad j(h)=hu
\lb{p}
\ee
$\R\,(\pt(hu)+\px(hu^2))$ and $\nu\px(h\px u)$ in
\rf{SSM1}--\rf{SSM2}, with $\R,\,\nu\ge 0$ denoting the Reynolds number and
viscosity parameter, represent inertial and Trouton viscosity terms,
respectively.

Additionally, the weak-slip and the strong-slip models 
contain as limiting cases three further lubrication models. 
One of them is the \emph{no-slip model}, which is obtained  setting $b=0$ in the weak-slip model:
\be
\pt h= - \px \Big(h^3\px \left(\pxx h-\Pi_\eps(h)\right)\Big) \,.
\label{NSM}
\ee
The second one is obtained from the 
strong-slip model in the limit $\beta\rightarrow\infty$ and describes the  
dynamics of suspended or falling free films:
\begin{subequations}
\label{FFM}
\begin{align}
\R\,(\eps\pt(hu)+\px(hu^2))&=\nu\px(h\px u)+h\px(\sigma\pxx h-\Pi_\eps(h))
\label{FFM1}\\
\eps\pt h &= -\,\px\left(h u\right),
\label{FFM2}
\end{align}
\end{subequations}
For the third limiting case  
the slip-length parameter $\beta_I$ is of order of magnitude 
lying in between those that lead to the weak 
and the strong-slip model, i.e. $b\ll \beta_I\ll\beta$. The corresponding 
\emph{intermediate-slip} model is given by 
\be
\pt h= - \px \Big(h^2\px \left(\pxx h-\Pi_\eps(h)\right)\Big) \,.
\label{ISM}
\ee
It can be obtained by rescaling time in \rf{WSM} by $b$ and 
letting $b\rightarrow\infty$ or by rescaling time and the horizontal velocity 
by $\beta$ in \rf{SSM1}--\rf{SSM2} and taking the limit $\beta\rightarrow 0$. Existence of weak solutions to \rf{SSM1}--\rf{SSM2} and \rf{FFM1}--\rf{FFM2} and rigorous convergence of the former ones to the classical solutions of \rf{ISM} as $\beta\go 0$ was shown recently in~\cite{KLN11}. 

As in~\cite{KLN11} we consider systems \rf{SSM1}--\rf{SSM2} and \rf{FFM1}--\rf{FFM2}  on a bounded interval $(0,L)$ with the boundary conditions
\be
u=0,\myand \px h=0 \myat x=0,\,L,
\lb{SSM BC}
\ee
whereas equations \rf{WSM},\rf{NSM} and \rf{ISM} with 
\be
\pxxx h=0,\myand \px h=0 \myat x=\pm L.
\label{WSM BC}
\ee
Both \rf{SSM BC} and \rf{WSM BC} incorporate zero flux at the boundary and as a consequence imply the
conservation of mass law
\bes
\frac{1}{L}\int_{0}^{L}h(x,t)\,dx=\const,\ \forall t>0.
\ees

Within the context of thin liquid films 
one of the first studies of the coarsening dynamics can be found in~\citet{GW03,GW04}. The authors considered the one-dimensional
no-slip lubrication model \rf{NSM} with \rf{WSM BC}. They confirmed numerically 
existence of the two coarsening driven mechanisms, namely collision and collapse. 
One of the typical problems considered in~\cite{GW03,GW04} was the calculation of the coarsening rates, i.e. how fast the number of
droplets decreases due to coarsening in time. Often in order to identify the characteristic
dependence for coarsening rates one needs to model very large arrays of
droplets (around $10^4$). But due to the presence of the ultrathin-layer of order $\eps$
between droplets the problem of numerical solution for any
lubrication equation becomes very stiff in time and demands high space resolution
as the number of droplets increases. 
Therefore, there exists a need for further reduction of lubrication models to more simple, possibly finite-dimensional ones.

Basing on the observation that solutions of lubrication equations describing
coarsening dynamics stay in time very close to a perturbed finite
combination of quasistationary droplets and can be therefore parameterized by a
finite number of parameters, namely positions and pressures of drops, in~\cite{GW03,GW04} for the first time
a reduced ODE model describing evolution of the latter ones on the slow time scale
was derived from the lubrication equation \rf{NSM}. Using this reduced model the authors
derived also the corresponding coarsening law in the form
\be
n(t)\sim t^{-2/5},
\lb{NCL}
\ee
where $n(t)$ denotes the number of droplets remaining at time $t$.
Later, analogous reduced ODE models from lubrication equation \rf{WSM} with a general mobility $M(h)=h^q$, $q>0$ in one and two dimensional case
were derived and analyzed in~\cite{ORS07,Gl07}.
An step to a rigorous justification of these models basing on a center manifold approach was made recently in~\cite{KRW11}.
For the case $M(h)=h$ the coarsening law \rf{NCL} was justified rigorously in~\cite{ORS06} using the gradient
flow structure of the corresponding lubrication equation. 
The work of~\cite{ORS07,Gl07} concerns migration of droplet. There it was shown that
the direction of the migration of droplets governed by  \rf{WSM} with a
general mobility $M(h)=h^q$, $q>0$ is opposite to the mass flux applied to
them. Moreover, for $g\le 2$ the driving coarsening mechanism is collapse of
droplets that is due to the mass diffusion in the ultra-thin layer between droplets and similar to Ostwald ripening in
binary alloys, see~\cite{BX94,BX95,SW00}. Note, that also in the no-slip case $q=3$,
i.e. one described by \rf{NSM}, as was shown in~\cite{GW04} the coarsening
rates even for the systems coarsening solely due collisions obey the law \rf{NCL}. 

Recently, in~\citet{KW09} it was shown that the coarsening dynamics of
quasistationary droplets governed by \rf{SSM1}--\rf{SSM2} with
sufficiently small $\mathrm{Re}$ number is driven also by collapse and collision. 
There reduced ODE models analogous to that one of~\cite{GW03,GW04} were
derived for system \rf{SSM1}--\rf{SSM2} and its limiting case
\rf{FFM1}--\rf{FFM2} as well.
In contrast, to the case of \rf{WSM} it was found there that the coefficients
of the strong-slip reduced ODE model depend explicitly on the slip length $\beta$.
In particular, there exists a critical length $\beta_{cr}=O(\eps)$ such that the migration of droplets proceeds in the direction of the applied mass flux
for $\beta>\beta_{cr}$ and opposite to it for $\beta<\beta_{cr}$. Moreover, it was shown that for moderate and large $\beta$ the driving coarsening mechanism   
switches from collapse to collision of droplets. Basing, on these observations
it was conjectured and shown numerically in~\cite{KW09} that the coarsening rates for 
systems \rf{SSM1}--\rf{SSM2} and \rf{FFM1}--\rf{FFM2} can be remarkably different from ones for \rf{NCL}.

In this study we continue the research initiated in~\cite{KW09}. 
Our aim here is to derive explicit coarsening laws for the dynamics of droplets in the strong-slip and free film regimes,
i.e. governed by lubrication system \rf{SSM1}--\rf{SSM2} and its limiting
case \rf{FFM1}--\rf{FFM2}. The missing point  in~\cite{KW09} was a derivation 
of flux representation between interacting droplets for moderate and large
slip lengths $\beta$ which was important for closure of the derived there reduced
ODE models. Therefore, inspired with the matched asymptotics technique applied
in~\citet{Gl07} to the lubrication equation \rf{WSM}, we present in section 2 a new closed form derivation
of reduced ODE models for \rf{SSM1}--\rf{SSM2} and \rf{FFM1}--\rf{FFM2} that
incorporate the explicit flux representation for all $0<\beta\le\infty$.

In section 3 we concentrate on the reduced ODE model corresponding to
\rf{FFM1}--\rf{FFM2}, i.e. on the regime of free films characterized
by the infinite slip length $\beta=\infty$.
In this case migration and subsequent collisions of droplets dominate
completely collapse component of the coarsening dynamics.
Therefore, we look only at the migration subsystem of the derived reduced  ODE
model such that droplet pressures are kept constant during evolution of
droplets and updated only after each subsequent collision event. We observe further that for a special initial data this migration subsystem 
can be solved explicitly while its solution represents subsequent collisions
of $N-1$ droplets with the largest last one. Therefore, we call it as an
exactly solvable collision/absorption model. 
It turns out that the coarsening law for this model depends only on the
initial distribution of the distances between droplets 
and can be derived analytically. Finally, we derive the continuous counterpart
of the coarsening law proceeding to the limit $N\go\infty$. 

In section 4 we consider several examples of initial distributions of distances between droplets and show that the corresponding coarsening rates
depend only on the distribution decay at infinity. 
Moreover, for an explicit family of distributions decaying as $1/x^{1+\alpha}$
with $\alpha>0$ we show existence of a threshold at $\alpha=1$ at which the
coarsening rates switch from algebraic to exponential ones. 

In section 5.1 we justify the derived hierarchy of the reduced models by numerical
comparison of their solutions to ones of the initial PDE system \rf{SSM1}--\rf{SSM2} 
and its limiting cases \rf{ISM} and \rf{SSM1}--\rf{SSM2}. We observe that the deviation between them stays
$O(\eps)$ uniformly in time.
Besides we compare solutions of the collision/absorption model from section 3
with those of the full reduced ODE system for the case $\beta=\infty$. Finally, in section
5.2 we check numerically 
the derived coarsening law for the collision/absorption model in the case of
finite $N$ and its continuous counterpart.

\section{Derivation of reduced ODE models.}\label{S.apriori}
We consider a solution to \rf{SSM1}--\rf{SSM2} which stays close in time to a union of $N+1$ droplets, which precise
characterization to be described below. Similar to the derivation of reduced coarsening models for the classical thin film equation in~\cite{Gl07}
we distinguish three regions in our matched asymptotic analysis.
\begin{itemize}
\item {\it Droplet core (DC) region:} This region corresponds to droplets and is composed of the union of disjoint intervals $(X_i(t)-R_i(t),X_i(t)+R_i(t))$,
so that $X_i(t)$ and $R_i(t)$ are the center and the radius of the $i$-th droplet, $i=0,...,N$.
The dynamical points $X_i(t)\pm R_i(t)$ are called contact line points and are defined through the relation
\be
h(X_i\pm R_i)=\eps H^*,
\lb{CLD}
\ee
where $H^*$ is the global maximum of $U'(H)$ with function $U(H)$ defined in \rf{U(h)}. We expand
\be
R_i=R_{i,0}+\eps R_{i,1}+..., X_i=X_{i,0}+\eps X_{i,1}+...
\lb{RXE}
\ee
and denote
\bes
\dot{R}=\frac{dR}{dt},\ \ \dot{X}=\frac{dX}{dt}.
\ees
\item {\it Contact line (CL) region} is a microscopic internal layer around the contact line points where $h$ and $x$ scale like $\eps$. Here we employ the moving rescaled spatial coordinate
\be
z=\frac{R(t)-|x-X(t)|}{\eps},\ \ \px=\frac{1}{\eps}\pz,\ \ \pt=\partial_\tau-\frac{1}{\eps}\pz(\dot{R}_0\mp\dot{X}_0),
\ee
where in this section the sign $\mp$ corresponds to two CL regions around the points $X\mp R$, respectively.
Accordingly, by definition \rf{CLD} we have $h(z=0)=\eps H^*$.
\item {\it Precursor layer (PL) region} is the complement $(-L,L)\setminus\cup_i(X_i-R_i,X_i+R_i)$. In this region $h$ scales like $\eps$.
\end{itemize}
The main goal is to determine the evolution of $R_i(t)$ and $X_i(t)$.
To do so as in~\cite{Gl07} we propose self-consistent asymptotic expansions in each of three regions and connect them via matching conditions. 
Corrections to the leading order base solutions solve linear equations, and Fredholm-type solvability conditions will yield information about the dynamics.

Let us first consider the motion of the $i$-th droplet with $i\in 1,...,N-1$. For a time being we skip below the subscript $i$.
Let us start with the CL region. Here the solution to \rf{SSM1}--\rf{SSM2} is expanded as
\bes
h=\eps H_1+\eps^2 H_2+...,\ \ u=\eps U_1+\eps^2 U_2+...
\ees
We will also use the induced expansions
\bes
P=P_0+\eps P_1+...,\ \ J=\eps^2 J_2+\eps^3 J_3+...
\ees
for the pressure and flux functions defined in \rf{p}.
The corresponding leading order system in $\eps$ in this region is given by
\beas
&&\pz(\sigma\pzz H_1-U'(H_1))=0,\\
&&\pz H_1(\dot{R}_0\mp\dot{X}_0)=\mp\pz(H_1 U_1).
\eeas
Integrating the last system and using matching conditions to the DC and PL regions
\bea
\pz H_1\go 0,\ \pzz H_1\go 0, H_1\go 1\ \text{as}\ \ z\go-\infty\nonumber\\
\pzz H_1\go 0,\ H_1\go +\infty\ \text{as}\ \ z\go+\infty
\lb{MC1}
\eea
one obtains
\sbea
\lb{H1}
\frac{\sigma}{2}(\pz H_1)^2&=&U(H_1)-U(1),\\
U_1&=&-(1-\frac{1}{H_1})(\dot{R}_0\mp\dot{X}_0)\mp\frac{J_2(-\infty)}{H_1}
\lb{U1}
\seea
In particular,
\be
\displaystyle \lim_{z\go+\infty}U_1=-(\dot{R}_0\mp\dot{X}_0),\ \ \lim_{z\go-\infty}U_1=\mp J_2(-\infty).
\ee

Next, in the DC region we expand the solution as
\bes
h=h_0+\eps h_1+\eps^2h_2+...,\ \ u=u_0+\eps u_1+\eps^2u_2+...
\ees
and correspondingly pressure as
\bes
p=p_0+\eps p_1+...
\ees
In turn, the leading order system in this region is given by
\beas
\sigma h_0\px(\pxx h_0)-\frac{u_0}{\beta}=0,\\
-\px(h_0u_0)=0
\eeas
Integrating the system and using the matching condition 
\be
h_0(X\mp R)=0
\lb{h0}
\ee
one obtains its solution in the form
\be
h_0=\frac{1}{R\sqrt{12\sigma}}(R^2-(x-X(t))^2),\ \ u_0\equiv 0.
\lb{DC}
\ee
Correspondingly, the leading order pressure is given by 
\be
p_0\equiv\frac{1}{R\sqrt{3\sigma}}.
\lb{PR}
\ee

In the PL region we expand the solution as
\bes
h=\eps h_1+\eps^2h_2+...,\ \ u=\eps u_1+\eps^2u_2+...
\ees
and correspondingly pressure and flux as
\bes
p=p_0+\eps p_1+...,\ \ j=\eps^2 j_2+\eps^3 j_3+...
\ees
The leading order system in this region is given by
\beas
&&h_1\px(U'(h_1))=0,\\
&&\pt h_1=-\px(h_1u_1)
\eeas
Integrating the system and using the matching condition $h_1(X\mp R)=1$ one obtains
\be
h_1\equiv 1,\ \ u_1=j_2\equiv\const
\lb{Ms01}
\ee

For the next order corrections $h_2,u_2$ in PL region one has the system
\beas
&&h_1\px(U''(h_1)h_2)=-\frac{u_1}{\beta},\\
&&\pt h_2=-\px(h_1u_2+h_2u_1)
\eeas
From the first equation and \rf{Ms01} one obtains 
\be
\pxx h_2=\pxx p_0=0\ \ \text{and}\ \ j_2=-\beta\px p_0.
\lb{j2}
\ee

Proceeding further in the expansion in the CL region for the second order corrections $H_2,U_2$ one obtains the system
\bea
-\nu\pz(H_1\pz U_1)&=&H_1(\sigma\pzz H_2-U''(H_1)H_2),\nonumber\\
0&=&\pz H_2(\dot{R}_0-\dot{X}_0)-\pz(H_1U_2+H_2U_1)
\lb{FOS_CL}
\eea
Let us introduce a linear operator
\bes
\boldsymbol{\mathcal{L}}\left[\begin{array}{ccc} 
H\\
U\end{array}\right]=\left[\begin{array}{ccc} 
H_1\pz(\sigma\pzz  H_2-U''(H_1)H_2))\\
\pz( H_2(\dot{R}_0-\dot{X}_0-U_1)-H_1U_2)\end{array}\right]\,.
\ees
Formal adjoint operator to $\boldsymbol{\mathcal{L}}$ is given by
\bes
\boldsymbol{\mathcal{L}}^*\left[\begin{array}{ccc}
g\\
 v\end{array}\right]=\left[\begin{array}{ccc}
-\sigma\pzzz(H_1 g)+U''(H_1)\pz(H_1g)-\pz V(\dot{R}_0-\dot{X}_0-U_1)\\
\pz v H_1
\end{array}\right].
\ees
The kernel of it contains two linear independent functions
\be
\left[\begin{array}{ccc}
g_{1}\\
 v_{1}\end{array}\right]:=
\left[\begin{array}{ccc}
1\\
0\end{array}\right],\ \
\left[\begin{array}{ccc}
g_{2}\\
 v_{2}\end{array}\right]:=
\left[\begin{array}{ccc}
1/H_1\\
0\end{array}\right].
\ee
To derive necessary Fredholm-type solvability conditions for the system \rf{FOS_CL} 
we multiply the first equation in \rf{FOS_CL} by $g_2$ and integrate it on $(-\infty,+\infty)$ to obtain
\bes
P_0(+\infty)-P_0(-\infty)=\int_{-\infty}^{+\infty}\frac{\nu}{H_1}\pz(H_1\pz U_1),
\ees
where we have used that in the CL region 
\be
P_0=U''(H_1)H_2-\sigma\pzz H_2.
\lb{P_0}
\ee
Substituting in the previous expression \rf{U1} one obtains
\be
P_0(+\infty)-P_0(-\infty)=-\nu I(\mp J_2(-\infty)+\dot{R}_0\mp\dot{X}_0),
\lb{PD}
\ee
where a constant integral $I$ is given by
\be
I=\int_{-\infty}^{+\infty}\frac{1}{H_1}\pz\left(\frac{\pz H_1}{H_1}\right)\,dz=\frac{1}{35(3 + \sqrt{3})}.
\lb{I}
\ee
and can be effectively calculated from \rf{H1} (see Appendix A). Formula \rf{PD} is an analog of Gibbs-Thomson boundary condition and
shows that the pressure experiences a jump at the CL region. Note, that this is a first considerable difference between the coarsening dynamics driven by
\rf{SSM1}--\rf{SSM2} and \rf{WSM}. In contrast to \rf{PD} as was shown in~\cite{GW03} the pressure is constant through the CL region in the case of \rf{WSM}.

Next, multiplying the first equation in \rf{FOS_CL} by $g_3$ and integrating it on $(-\infty,+\infty)$ one obtains
\bes
0=\nu H_1\pz U_1\Big|_{-\infty}^{+\infty}+\int_{-\infty}^{+\infty}H_1\pz(\sigma\pzz H_2-U''(H_1)H_2).
\ees
Integrating further three times by parts and using \rf{U1}, \rf{P_0} one arrives at
\beas
0&=&-\nu\frac{\pz H_1}{H_1}(\mp J_2(-\infty)+\dot{R}_0\mp\dot{X}_0)\Big|_{-\infty}^{+\infty}-H_1P_0\Big|_{-\infty}^{+\infty}\\[2ex]
&-&\sigma\pz H_1\pz H_2\Big|_{-\infty}^{+\infty}+\sigma\pzz H_1H_2\Big|_{-\infty}^{+\infty},
\eeas
Using the matching condition \rf{MC1} and additionally 
\bea
&&\pz H_2\go \const\ \text{as}\ \ z\go-\infty,\nonumber\\
&&\pz H_1\go \px h_0,\ \pz H_2\sim\px h_1+\pxx h_0 z,\ H_1\sim h_1+\px h_0 z\ \text{as}\ \ z\go+\infty
\lb{MC2}
\eea
one arrives at
\bes
(H_1P_0)\Big|_{-\infty}^{+\infty}=\sigma\pz H_1(+\infty)\pz H_2(+\infty).
\ees
The last expression again using \rf{MC1} and  \rf{MC2} implies
\bea
\sigma\pxx h_0&=&-P(+\infty),\nonumber\\
\sigma(\px h_0\px h_1)\Big|_{X\mp R}&=&P(-\infty)-P(+\infty)h_1(X\mp R).
\lb{MC3}
\eea
Note, that the first relation in \rf{MC3} is consistent with already derived \rf{DC}--\rf{PR},
whereas the second one is new.

Finally, let us consider the system for the first order corrections $h_1,u_1$ in the DC region which has the form
\sbea
\lb{FOS_DC_u}
0&=&\nu\px(h_0\px u_1)+\sigma h_0\pxxx h_1-u_1/\beta,\\
\frac{\partial h_0}{\partial R}\dot{R}_0-\frac{\partial h_0}{\partial x}\dot{X}_0&=&-\px(h_0u_1) 
\lb{FOS_DC_h}
\seea
Let us introduce a linear operator
\bes
\boldsymbol{\mathcal{L}}\left[\begin{array}{ccc} 
h\\
u\end{array}\right]=\left[\begin{array}{ccc} 
\nu\px(h_0\px u_1)+\sigma h_0\pxxx h_1-u_1/\beta\\
-\px(h_0u_1)\end{array}\right]\,.
\ees
Formal adjoint operator to $\boldsymbol{\mathcal{L}}$ is given by
\bes
\boldsymbol{\mathcal{L}}^*\left[\begin{array}{ccc}
g\\
 v\end{array}\right]=\left[\begin{array}{ccc}
\nu\px(h_0\px g)-\frac{g}{\beta}+h_0\px v\\
-\sigma\pxxx(h_0g)
\end{array}\right].
\ees
The kernel of it contains two linear independent functions
\be
\left[\begin{array}{ccc}
g_{1}\\
 v_{1}\end{array}\right]:=
\left[\begin{array}{ccc}
0\\
1\end{array}\right],\ \ 
\left[\begin{array}{ccc}
g_{2}\\
 v_{2}\end{array}\right]:=
\left[\begin{array}{ccc}
1\\
\int_X^x\frac{d\tau}{\beta h_0}\end{array}\right].
\ee
To derive necessary Fredholm-type solvability conditions for the system \rf{FOS_DC_u}--\rf{FOS_DC_h} 
we multiply \rf{FOS_DC_h} by $v_1$, integrate it
and using the matching condition \rf{h0} obtain
\be
\dot{R}_0=0.
\lb{R0}
\ee
In turn, multiplying the right hand side of the second equation in \rf{FOS_CL}
by $v_2$ and integrating it on $(X-R,\,X+R)$ one obtains
\beas
0&=&-\dot{X}_0\int_{X-R}^{X+R}\frac{\partial h_0}{\partial x}v_2\,dx+\int_{X-R}^{X+R}\px(h_0u_1)v_2\,dx=h_0u_1v_2\Big|_{X-R}^{X+R}-\int_{X-R}^{X+R}h_0u_1\px v_2\,dx-\\
&-&\dot{X}_0\int_{X-R}^{X+R}\frac{\partial h_0}{\partial x}v_2\,dx=\dot{X}_0(h_0v_2)\Big|_{X-R}^{X+R}-\frac{2\dot{X}_0R}{\beta}-\frac{1}{\beta}\int_{X-R}^{X+R}u_1\,dx=\\
&=&-\frac{2\dot{X}_0R}{\beta}-\frac{1}{\beta}\int_{X-R}^{X+R}u_1\,dx.
\eeas
In the last equality we used that $h_0\sim O(R-|x-X|)$ and $v_2\sim\log(R-|x-X|)$ as $x\go X\mp R$.
Next, using \rf{FOS_DC_u} and integrating three times by parts one arrives at
\bea
\frac{2\dot{X}_0R}{\beta}&=&-\frac{1}{\beta}\int_{X-R}^{X+R}u_1\,dx=\int_{X-R}^{X+R}\nu\px(h_0\px u_1)+\sigma h_0\pxxx h_1\,dx=\nonumber\\
&=&\left[\nu h_0\px u_1+\sigma h_0\pxx h_1-\sigma\px h_0\px h_1+\sigma\pxx h_0 h_1\right]\Big|_{X-R}^{X+R}
\lb{Ms1}
\eea 
Let us note that from \rf{FOS_DC_u}--\rf{FOS_DC_h} and \rf{R0}, \rf{h0} it follows that 
\sbea
\lb{u1}
u_1\equiv\dot{X}_0,\\
\pxxx h_1=\frac{\dot{X}_0}{\beta\sigma h_0}
\lb{h1}
\seea
Hence, by \rf{u1} the first term in the square brackets in \rf{Ms1} vanishes. By \rf{h1} one has
$\pxx h_1\sim\log(R-|x-X|)$ as $x\go X\pm R$. Due to this and \rf{h0}--\rf{DC} the second term  in the square brackets in \rf{Ms1} also vanishes.
In turn, due to the matching condition 
\be
h_1(X\mp R)=H_1(0)
\lb{MC4}
\ee
and \rf{DC} the fourth fourth term  in the square brackets in \rf{Ms1} vanishes.
Therefore, relation \rf{Ms1} reduces to
\be
\frac{2\dot{X}_0R}{\beta}=-\left[\sigma\px h_0\px h_1\right]\Big|_{X-R}^{X+R}.
\lb{Ms2}
\ee

At this moment let us introduce back the droplet subscript $i=1,...,N-1$ and denote by $J_i$ the flux $j_2$ in the PL region between
$i-1$-th and $i$-th droplets. Combining \rf{j2} with \rf{PD} and \rf{R0} one obtains that $J_i$ is constant and satisfies
\be
J_i=\beta\frac{P_i-P_{i-1}-2\nu J_iI-\nu I(\dot{X}_{i,0}+\dot{X}_{i-1,0})}{d_i}.
\lb{Ms3}
\ee
In the last formula we introduced two more notations: the constant leading order pressure inside $i$-th droplet
\be
P_i=\frac{1}{R_i\sqrt{3\sigma}}.
\lb{P}
\ee
according to \rf{PR} and the distance between the neighboring DC regions
\bes
d_i=X_i-X_{i-1}-R_i-R_{i-1}.
\ees 
From \rf{Ms3} one obtains an explicit expression for $J_i$:
\be
J_i=\beta\frac{P_i-P_{i-1}-\nu I(\dot{X}_{i,0}+\dot{X}_{i-1,0})}{d_i+2\nu I\beta},\ \ i=1,...,N.
\lb{FR}
\ee
Next, from \rf{Ms2}, the matching conditions \rf{MC3}, \rf{MC4} and equations\rf{PD}, \rf{R0} one obtains
the leading order equation for $i$-th droplet position evolution
\bes
\dot{X}_{i,0}=-\frac{I\beta\nu}{2R_i+2I\beta\nu}(J_{i+1}+J_i).
\ees
Substituting in the last expression the flux representation \rf{FR}, denoting
\be
\widetilde{d}_i=\frac{d_i}{I\nu\beta},
\lb{Cn}
\ee
and using \rf{P} one obtains
\bea
\hspace{-0.3cm}\dot{X}_{i,0}&=&-\frac{P_i}{2/(\sqrt{3\sigma}\beta)+2I\nu P_i}\left(\frac{(P_{i+1}-P_i)-I\nu(\dot{X}_{i+1}+\dot{X}_i)}{\widetilde{d}_{i+1}+2}
+\frac{(P_i-P_{i-1})-I\nu(\dot{X}_i+\dot{X}_{i-1})}{\widetilde{d}_i+2}\right),\nonumber\\[2ex]
&&\text{for}\ \ i=1,...,N-1.
\lb{X0}
\eea
In turn by \rf{R0} and definition \rf{P} one has
\bes
\dot{P}_{i,0}=0,\ \ \text{for}\ \ i=0,...,N.
\ees
The derived ODE system describing the leading order in $\eps$ evolution of pressures and positions of $N+1$ droplets will be closed if we additionally prescribe
that the first and the last droplet do not move, i.e
\be
\dot{X}_{0,0}=\dot{X}_{N,0}=0,\ X_0=0,\,X_{N}=L. 
\lb{BC}
\ee
The condition \rf{BC} corresponds to the situation when one extends the array on $N+1$ droplets from interval $(0,L)$ to an infinite array on the whole real line $\mR$ by reflection around the points $x=0$ and $x=L$. It stays also in agreement with boundary conditions \rf{SSM BC} and \rf{WSM BC}

Let us point out that the evolutions of pressures is slower than one of positions and proceeds on the order $\eps$. One can potentially obtain it by going further
in the expansion of the solution to \rf{SSM1}--\rf{SSM2}, while an easier way is to derive it from the conservation of droplet volume as was done 
in~\cite{GW03,Gl07} for the case of equation \rf{WSM}. Namely, the volume of the $i$-th droplet $V_i$ is changing due to the difference of the fluxes in the surrounding it PL regions. Using \rf{DC} and \rf{P} one obtains
\bes
\dot{V}_{i,0}=\eps\frac{4A^3}{3P^3}\dot{P}_{i,1}=\eps^2(J_{i+1}-J_i),
\ees
Substituting in the last expression the flux representation \rf{FR} and denoting
\bes
C_i=\eps\frac{3P^3}{4A^3}
\ees one obtains
\bes
\eps\dot{P}_{i,1}=\frac{C_i}{I\nu}\left(\frac{P_{i+1}-P_i}{\widetilde{d}_{i+1}+2}-\frac{P_i-P_{i-1}}{\widetilde{d}_i+2}\right)-
C_i\left(\frac{\dot{X}_{i+1}-\dot{X}_i}{\widetilde{d}_{i+1}+2}-\frac{\dot{X}_i-\dot{X}_{i-1}}{\widetilde{d}_i+2}\right),\ i=1,...,N-1. 
\ees
Finally, combining the last expression with \rf{X0} and \rf{BC} the closed ODE system for the leading order evolution of positions and pressures in the array of $N+1$ droplets takes the following form:
\bea
\dot{X}_i&=&-\frac{P_i}{2/(\sqrt{3\sigma}\beta)+2I\nu P_i}\left(\frac{(P_{i+1}-P_i)-I\nu(\dot{X}_{i+1}+\dot{X}_i)}{\widetilde{d}_{i+1}+2}
+\frac{(P_i-P_{i-1})-I\nu(\dot{X}_i+\dot{X}_{i-1})}{\widetilde{d}_i+2}\right),\nonumber\\[2ex]
\dot{P}_i&=&\frac{C_i}{I\nu}\left(\frac{P_{i+1}-P_i}{\widetilde{d}_{i+1}+2}-\frac{P_i-P_{i-1}}{\widetilde{d}_i+2}\right)-
C_i\left(\frac{\dot{X}_{i+1}-\dot{X}_i}{\widetilde{d}_{i+1}+2}-\frac{\dot{X}_i-\dot{X}_{i-1}}{\widetilde{d}_i+2}\right),\ i=1,...,N-1;
\lb{ODE}
\eea
and
\bea
\dot{P}_1+2C_1\frac{\dot{X}_2}{\widetilde{d}_1+2}&=&2\frac{C_1}{I\nu}\frac{P_2-P_1}{\widetilde{d}_1+2},\ \ \dot{X}_1=0,\nonumber\\
\dot{P}_N-2C_N\frac{\dot{X}_2}{\widetilde{d}_{N-1}+2}&=&-2\frac{C_N}{I\nu}\frac{P_N-P_{N-1}}{\widetilde{d}_{N-1}+2},\ \ \dot{X}_N=0.
\lb{PB}
\eea
Let us consider certain limiting cases.
In the case $\beta\go\infty$ the limiting system for evolution of pressures and positions has the form
\bea
\dot{P}_i+C_i(\dot{X}_{i+1}-\dot{X}_{i-1})&=&\frac{C_i}{I\nu}(P_{i+1}-2P_i+P_{i-1}),\nonumber\\
\dot{X}_{i+1}-2\dot{X}_i+\dot{X}_{i-1}&=&\frac{P_{i+1}-P_{i-1}}{\nu I},\ \ \text{for}\ \ i=1,...,N-1;
\lb{ODEinf}
\eea
and 
\bea
\dot{P}_1+C_1\dot{X}_2&=&\frac{C_1}{I\nu}(P_2-P_1),\ \ \dot{X}_1=0,\nonumber\\
\dot{P}_N-C_N\dot{X}_{N-1}&=&-\frac{C_N}{I\nu}(P_N-P_{N-1}),\ \ \dot{X}_N=0.
\lb{PBinf}
\eea
Next, rescaling the time by $\beta\nu$ and proceeding to the limit $\beta\go 0$  the limiting system for evolution of pressures and positions  
takes the form
\bea
\dot{P}_i&=&C_i\left(\frac{P_{i+1}-P_i}{d_{i+1}}-\frac{P_i-P_{i-1}}{d_i}\right),\nonumber\\
\dot{X}_i&=&-\frac{P_i\sqrt{3\sigma}I}{2}\left(\frac{P_{i+1}-P_i}{d_{i+1}}+\frac{P_i-P_{i-1}}{d_i}\right),\ \ \text{for}\ \ i=1,...,N-1;
\lb{ODE0}
\eea
and 
\bea
\dot{P}_1&=&2C_1\frac{P_2-P_1}{d_1},\ \ \dot{X}_1=0,\nonumber\\
\dot{P}_N&=&-2C_N\frac{P_N-P_{N-1}}{d_N},\ \ \dot{X}_N=0.
\lb{PB0}
\eea
Note, that the last system coincides with one derived in~\cite{Gl07} for the intermediate-slip equation \rf{ISM} in the one-dimensional case.
This stays in agreement with the fact that \rf{ISM} is the limiting case of \rf{SSM1}--\rf{SSM2} as $\beta\go 0$ as was shown in~\cite{MWW06,KLN11}.
Finally, note that after time rescaling by $\beta\nu$ taking limits $\nu\go \infty$ or  $\nu\go 0$ results again in \rf{ODEinf}--\rf{PBinf} and 
\rf{ODE0}--\rf{PB0}, respectively. This is also
natural, because \rf{FFM1}--\rf{FFM2} and \rf{ISM} are the limiting cases  of
\rf{SSM1}--\rf{SSM2} as well as  $\nu\go \infty$ or  $\nu\go 0$, respectively. 
 
Let us summarize the algorithm for simulation of coarsening dynamics in large arrays of droplets using the derived reduced ODE models.
Starting with an array of $N+1$ droplets 
after each subsequent coarsening event (i.e a collapse of one droplet or collision of two droplets) one can model the coarsening process
further by reducing the dimension of the model by two and solving the reduced ODE model with the updated initial data. 
Practically, as in~\cite{GW04} we say that a collapse event occurs at a moment when pressure of one droplet increases a certain threshold, namely when 
\be
P>0.5P_{max}(\eps),\ \ \text{with}\ \ P_{max}(\eps):=\frac{27}{256\eps}.
\lb{Pmax}
\ee
Then we take the final pressures and positions for remaining droplets from the previous run of the reduced ODE model
as initial conditions for the next one. In the case of collision in~\citet{GW04} was suggested that coarsening event occurs when
the distance between two colliding $i$-th and $i+1$-th droplets becomes
smaller then a certain threshold $\delta=O(\eps)$, i.e. when
\be
d_i\le \delta,
\label{Coldist}
\ee
After the collision we calculate the position and the pressure for the new formed  droplet by formulas
\begin{align}
X_{i,new}&=1/2(X_{i+1}-R_{i+1}+X_i-R_i),\nonumber\\
P_{i,new}&=\left(\frac{1}{P_i^2}+\frac{1}{P_{i+1}^2}\right)^{-1/2}.
\label{ClnRecal}
\end{align}
The last formula for $P_{i,new}$ is based on the observation that mass of the new droplet
is to the leading order in $\eps$ given by the sum of the masses of the collided droplets (see~\cite{GW04}).
In section 5.1 we compare solutions of the derived reduced ODE model \rf{ODE}--\rf{PB} with those of the initial PDE system \rf{SSM1}--\rf{SSM2} and show
that the former ones provide high accuracy $O(\eps)$ also after subsequent coarsening events.

\section{An exactly solvable collisions/absorption model} 

Let us consider the limiting case of infinite slip length $\beta=\infty$, namely the ODE system \rf{ODEinf}--\rf{PBinf} describing coarsening in free films. 
As pressure evolution proceeds on a slower time scale then that one of positions as $\eps\go 0$ let us
consider only migration of droplets. Namely, we investigate the zero order system 
\bea
\dot{X}_0&=&\dot{X}_N=\dot{P}_i=0,\ \ \text{for}\ \ i=0,...,N,\ \;\nonumber\\
\dot{X}_{i+1}-2\dot{X}_i+\dot{X}_{i-1}&=&\frac{P_{i+1}-P_{i-1}}{\nu I},\ \ \text{for}\ \ i=2,...,N-1;
\lb{RS}
\eea
As will be justified numerically in section 5.1 for given $\eps,\,T>0$ one can choose initial data with sufficiently small $P_i(0)\ll 1,\ i=0,1,...,N$
such that the difference between solutions to \rf{RS} and \rf{ODEinf}--\rf{PBinf} stays uniformly $O(\eps)$ for all times $t\in (0,T]$. 
Note that for such initial data there is no other constraint on the location of $X_i(0)$ other than that $d_i(0)$ should be larger then the collision threshold 
introduced in \rf{Coldist}.

Moreover, for certain initial data one can solve \rf{RS} explicitly. Indeed, if
\be
P_i(0)=p,\ \text{for}\ i=0,1,...,N-1\ \text{and}\ P_N(0)=\bar{p}\ \text{with}\ 1\gg p>\bar{p}
\lb{ID}
\ee
then the solution to \rf{RS} is given by
\be
X_i(t)=X_i(0)+\frac{Bi}{N}t,\ \text{for}\ i=1,...,N-1;\ X_0=0,\,X_N=L,\ \text{where}\ B=\frac{p-\bar{p}}{\nu I}.
\lb{XSS}
\ee 
In this and the next sessions for convenience reasons we redefine the notation for the distances between droplets from the previous section as
\bes
d_i(t)=X_i(t)-X_{i-1}(t)\ \text{for}\ i=1,...,N-1\ \text{and}\ d_N(t)=L-X_{N-1}(t)-R_N(t)-R_{N-1}(t)
\ees
and call usually $d_i(t)$ as the distance of the $i$-th droplet at time $t$.
Using this notation one can rewrite the solution \rf{XSS} in the following form
\bea
d_i(t)&=&d_i(0)+\frac{B}{N}t,\ \ \text{for}\ \ i=1,...,N-1\nonumber\\
d_N(t)&=&d_N(0)-\frac{B(N-1)}{N}t\ \text{for}\ t\in (0,T_c),\ \text{where}\ T_c=\frac{d_NN}{(N-1)B}.
\lb{RSS} 
\eea
Note, that $T_c$ denotes the time proceeding until $(N-1)$-th droplet collides with the largest last one.
Iterating \rf{RSS} one observes that first $N-1$ droplets collide one after another with the last one like rings in the famous rubber band toy. Due to \rf{XSS} all
droplets except the first and the last ones move to the right.  The last droplet  consequently absorbs the neighbor droplet, while the distance
between them is uniformly distributed between the remaining droplets. Therefore, the distances of the remaining 
droplets at the collision time $T_c$ are given by
\be
d_i(T_c)=d_i(0)+d_N(0)/(N-1),\ \ i=1,...,N-1.
\lb{hu}
\ee
Writing the solution to \rf{RS} in the form \rf{RSS} is convenient because one
can substitute $d_i(T_c)$ as the initial distances for the modelling of the next collision event.

Note, that due to \rf{hu} distance monotonicity is preserved in time for solutions \rf{RSS}, i.e. if $d_l(0)>d_m(0)$ for some $l,m\in 1,...,N$ then
$d_l(t)>d_m(t)$ for all times $t>0$. This property allows us basing only on a given initial distribution of the distances in the array of droplets
to derive the coarsening laws analytically for solutions to \rf{RS} considered with \rf{ID} and additional assumption
\be
1>>p>>P_N.
\lb{ID2}
\ee
This assumption prescribes that the last droplet is much larger then others 
and allows us to simplify further the dynamics by assuming that its pressure $P_N$ remains constant in time.
In turn, this implies that the coarsening dynamics in this case depends solely on the evolution of droplet positions without change of their pressures after subsequent collisions.

Indeed, let initial distances be prescribed by $k\in\N$ families such that there are $i_m$ distances in $m$-th family ($1\le m\le k$), 
all of them are equal to $d_m$ and
\be
d_1\ge d_2\ge...\ge d_k,\quad i_1+i_2+...+i_k=N
\lb{hi}
\ee
holds. Additionally, let us allocate these $k$ families in the initial configuration so that the distances between droplets non increase coming from 
the first to the last droplet. Then due to the distance monotonicity property this ordering will be preserved in time, i.e. first the members of the family $k$ will be absorbed by the last droplet, then those of $k-1$-th one and etc. Moreover, the distances in each family will stay equal for all $t>0$. This implies that for initial data  satisfying \rf{ID}, \rf{ID2} and \rf{hi} all collision times are uniquely determined having given  $k$ and the set $\{d_m,i_m\},\ m=1,...,k$. Therefore, using the explicit solution \rf{RSS} holding between subsequent collision events the coarsening law can be derived analytically by a recursive procedure.

Indeed, let us fix an index $1\le m\le k$ and look at the moment when all families with the indexes $m+1,...,k$ and also $l-1$ members of the $m$-th family
have been absorbed for some given $1\le l\le i_m$. Let us calculate the time $t(n)$ needed for the absorption of the $l$-th member with $n$ denoting
the remaining number of droplets after the latter event. Using \rf{RSS} one can easily calculate by recursion that
\be
\displaystyle t(n)=\frac{n+1}{nB}\left(\widetilde{d}_m+\sum_{r=1}^{l-1}\frac{t(n+r)B}{n+r+1}\right)=\frac{(n+l)\widetilde{d}_m}{nB},
\lb{Ms03}
\ee
where by $\widetilde{d}_m$ we denote the distance in the $m$-th family at the time when $(m+1)$-th family has been absorbed. 
From \rf{Ms03} one can obtain the total time needed for the $m$-th family to be absorbed $T_m$ in the form
\be
\displaystyle T_m=\frac{\widetilde{d}_m}{B}\left(N-\sum_{p=m+1}^{k}i_p\right)\sum_{r=1}^{i_p}\frac{1}{N-\sum_{p=m+1}^{k}i_p-r}.
\lb{Tm}
\ee    
In turn, using again \rf{RSS} one recursively finds
\bes
\displaystyle \widetilde{d}_m=\sum_{p=m+1}^{k}\frac{ \widetilde{d}_pi_p}{N-\sum_{p'=p}^{k}i_p}=\frac{\sum_{p=m+1}^{k}d_pi_p}{N-\sum_{p'=p}^{k}i_p}+d_m
\ees
Substituting the last expression in \rf{Tm} one obtains
\bes
\displaystyle T_m=\frac{1}{B}\left(Nd_m+\sum_{p=m}^{k}(d_p-d_m)i_p\right)\sum_{r=1}^{i_m}\frac{1}{N-\sum_{p=m+1}^{k}i_p-r}
\ees
Therefore, the total time needed for all families up to $m$-th to be absorbed is given by 
\be
\displaystyle T(d_m)=\sum_{p=m}^{k}\frac{1}{B}\left(Nd_p+\sum_{p'=p}^{k}(d_p'-d_p)i_p'\right)\sum_{r=1}^{i_p}\frac{1}{N-\sum_{p'=p+1}^{k}i_p'-r}
\lb{DL}
\ee

Let us now derive the continuum version for the discrete coarsening law in \rf{DL} proceeding to the limit $N\go\infty$ and $k\go\infty$.
Suppose we are given a probability density function $f(d)$ on $(0,+\infty)$, i.e.
\bes
\int_0^{+\infty} f(x)\,dx=1,\ \ f(d)\ge 0\ \text{and}\ f(d)=0\ \text{if}\ d\le 0.
\ees
Defining $d_i=i\Delta d$ for $i\in\N\cup\{0\}$ and a fixed $\Delta d\ll 1$ we approximate $f(d)$ by a piece-wise constant function $f_a(d)$ as follows.
\bes
f_a(d)=f(d_{i+1})\ \text{for}\ d\in [d_i,\,d_{i+1})\ \text{and}\ f_a(d)=0\ \text{if}\ d\le 0. 
\ees
Accordingly to this approximation suppose we are given an array of $N+1$ droplets with $N\gg 1$ such that the number of droplets with the distances lying in the interval
$[d_i,\,d_{i+1})$ is equal to $[Nf(d_{i+1})\Delta d]$. As before we suppose that droplets are allocated so
that the distances between them are non increasing and \rf{ID}, \rf{ID2} hold. Then using \rf{DL} one obtains
\beas
\hspace{-1.3cm}T_{\Delta d, N}(d_m)&=& 
\sum_{p=0}^{m}\frac{1}{B}\left(Nd_p+\sum_{p'=p}^{k}(d_p'-d_p)Nf(d_p')\Delta d\right)\sum_{r=1}^{Nf(d_p)\Delta d}\frac{1}{N-\sum_{p'=p+1}^{k}Nf(d_p')\Delta d-r}+O(\Delta d, 1/N)
\\
&=&\sum_{p=0}^{m}\frac{1}{B}\left(d_p+\sum_{p'=p}^{k}(d_p'-d_p)f(d_p')\Delta d\right)\sum_{r=1}^{Nf(d_p)\Delta d}\frac{1}{1-\sum_{p'=p+1}^{k}f(d_p')\Delta d-r/N}+O(\Delta d, 1/N)
\\ &=&
\sum_{p=0}^{m}\frac{N}{B}\left(d_p+\sum_{p'=p}^{k}(d_p'-d_p)f(d_p')\Delta d\right)\sum_{s=1/N}^{f(d_p)\Delta d}\frac{\Delta s}{1-\sum_{p'=p+1}^{k}f(d_p')\Delta d-s}+O(\Delta d, 1/N)
\eeas
Taking the limit $N\go+\infty$ in the last expression and introducing
\bes
\displaystyle T_{\Delta d}(d)=\lim_{N\go+\infty}\frac{T(d)_{\Delta d, N}}{N}
\ees
one obtains
\be
T_{\Delta d}(d)=\sum_{p=0}^{m}\frac{N}{B}\left(d_p+\sum_{p'=p}^{k}(d_p'-d_p)f(d_p')\Delta d\right)\int_{s=0}^{f(d_p)\Delta d}\frac{ds}{1-\sum_{p'=p+1}^{k}f(d_p')\Delta d-s}+O(\Delta d).
\lb{Ms6}
\ee
Applying the Taylor expansion to the last integral in \rf{Ms6} one finds
\bes
\int_{s=0}^{f(d_p)\Delta d}\frac{ds}{1-\sum_{p'=p+1}^{k}f(d_p')\Delta d-s}=\frac{f(d_p)\Delta d}{1-\sum_{p'=p}^{k}f(d_p')\Delta d}+O(\Delta d^2)
\ees
Inserting this into \rf{Ms6} one obtains
\bes
T_{\Delta d}(d)=\frac{1}{B}\sum_{p=0}^{m}\left(d_p+\sum_{p'=p}^{k}(d_p'-d_p)f(d_p')\Delta d\right)\frac{f(d_p)\Delta d}{1-\sum_{p'=p}^{k}f(d_p')\Delta d}+O(\Delta d).
\ees
Finally, taking the limit $\Delta d\go 0$ and introducing
\bes
\displaystyle T(d)=\lim_{\Delta d\go 0}\frac{T_{\Delta d}(d)}{N}
\ees
one arrives at
\be
T(d)=\frac{1}{B}\int_{0}^{d}\left(x+\int_{0}^{x}(y-x)f(y)\,dy\right)\frac{f(x)}{1-\int_{0}^{x}f(y)\,dy}\,dx
\lb{Ms7}
\ee
Introducing function $n(d)$ as the relative number of droplets with initial distances larger or equal $d$, i.e. as
\be
n(d)=1-\int_0^d f(x)\,dx
\lb{nh}
\ee
one obtains from \rf{Ms7} that
\bea
T(d)&=&\frac{1}{B}\int_{0}^{d}\left(-x+\int_{0}^{x}(y-x)n'(y)\,dy\right)\frac{n'(x)}{1-\int_{0}^{x}n(y)\,dy}\,dx=\nonumber\\
&=&\frac{1}{B}\int_{0}^{d}n(x)\ln\left[\frac{n(x)}{n(d)}\right]\,dx.
\lb{CRL}
\eea
The last expression provides an exact coarsening law, i.e. it tells what time $T(d)$ is needed 
until all droplets having initially distances smaller then $d$ are absorbed by the last large droplet.

In appendix B, we show that the discrete coarsening law \rf{DL} can be recovered back from \rf{CRL} if the initial distribution $f(x)$ 
has the form
\be
f(d)=\sum_{m=1}^ki_m'\delta(d-d_m),
\lb{DD}
\ee
i.e. if initial distance distribution is represented by $k\in\N$ families as in \rf{hi} while the number of droplets $N\go\infty$.
Moreover, we justify the connection between \rf{CRL}, \rf{DL} and the starting ODE system \rf{RS} as well numerically in the section 5.2. 

\section{Examples of coarsening rates} 

{\bf a)} We consider an explicit family of initial distributions $f(x)$ and show that depending on their decay as
$x\go+\infty$ the coarsening rates reproduce all possible algebraic decays. Moreover, there is a certain threshold after which the decay becomes exponential.
Namely, let us consider a family 
\be
f(x)=\frac{1}{x^{1+\alpha}}\Big/\int_{A}^{+\infty}\frac{dx}{x^{1+\alpha}}=\frac{A^\alpha}{x^{1+\alpha}}\ \ \text{with}\ \ \alpha,\,A>0.
\lb{PDF}
\ee
From \rf{nh} it follows that 
\be
n(x)=\left(\frac{A}{x}\right)^\alpha.
\lb{nha}
\ee  
Substituting this in \rf{CRL} one obtains
\sbea
\lb{alpha}
T(d)&=&\frac{\alpha A}{B(\alpha-1)}\left(\frac{1}{\alpha-1}\left[\left(\frac{d}{A}\right)^{1-\alpha}-1\right]+\alpha\ln\left[\frac{d}{A}\right]\right)\ \text{if}\ \alpha\ne 1,\\
T(d)&=&\frac{A}{B}\left(\left[\ln\left(\frac{d}{A}\right)\right]^2/2+\ln\left(\frac{d}{A}\right)\right)\ \text{if}\ \alpha=1.
\lb{alpha1}
\seea
Combining \rf{alpha1} and \rf{nha} one obtains the exact coarsening law for the case $\alpha=1$
\be
n(t)=\exp\left[1-\sqrt{1+2Bt/A}\right]
\lb{CR1}
\ee
In the case $\alpha\ne 1$ one obtains from \rf{alpha} and \rf{nha} 
\bes
T(n)=\frac{\alpha A}{B(\alpha-1)}\left(\frac{1}{1-\alpha}\left[n^\frac{\alpha-1}{\alpha}-1\right]+\ln(n)\right).
\ees
For the latter exact law one obtains the following asymptotics
\be
n(t)\sim
\left\{
\begin{aligned}
\displaystyle &\left(\frac{tB(\alpha-1)^2}{\alpha A}\right)^\frac{\alpha}{\alpha-1},\hspace*{2cm}\text{if}\ \ \alpha<1\\[2ex]
\displaystyle &\exp\left\{-\frac{tB(\alpha-1)}{\alpha A}\right\},\hspace*{2cm}\text{if}\ \ \alpha>1
\end{aligned}\right.\ \ \text{as}\ \ t\go\infty.
\lb{CR}
\ee
Therefore, from \rf{CR1}--\rf{CR} one finds out that for $0<\alpha<1$ the coarsening rates are algebraic at least for large times, while at $\alpha=1$ they become exponential and stay so for $\alpha\in (1,\,+\infty)$.

{\bf b)} Consider $f(x)=\exp(-x)$. Substituting it in \rf{nh} and \rf{CRL}, consequently, one obtains the exact law
\bes
T(n)=\frac{1}{B}(n-1-\ln(n)).
\ees
Thus, in this case the following asymptotics holds
\be
n(t)\sim\exp(-Bt)\ \ \text{as}\ \ t\go\infty.
\lb{CRexp}
\ee

{\bf c)} Consider a Gaussian distribution $f(x)=2/\sqrt{\pi}\exp(-x^2)$. In this case by \rf{nh} one has $n(x)=\mathrm{erfc}(x)$. Substituting it in
\rf{CRL} one obtains
\beas
T(d)&=&\frac{1}{B}\int_0^d\left[\int_0^xn(y)\,dy\right]\frac{n'(x)}{n(x)}\,dx=\\
&=&\frac{1}{B\sqrt{\pi}}\int_0^d\left(\frac{(1-\exp(-x^2))\exp(-x^2)}{\int_x^{+\infty}\exp(-x^2)\,dt}-2x\exp(-x^2)\right)\,dx\\
&=&\frac{1}{B\sqrt{\pi}}\left(-C+O\left(\exp(-d^2)\right)+\int_0^d\frac{\exp(-x^2)}{\int_x^{+\infty}\exp(-x^2)\,dt}\,dx\right)\\
&=&\frac{1}{B\sqrt{\pi}}(-C-\ln n(d))+O(\exp(-d^2)),
\eeas
where constant $C\approx 0.74$.
Therefore the following asymptotics holds
\be
n(t)\sim\exp\{-C-B\sqrt{\pi}t\}\ \ \text{as}\ \ t\go\infty.
\lb{CRgauss}
\ee
and the coarsening rates show an exponential decay as in the example b).
\begin{figure}
\includegraphics*[width=0.5\textwidth]{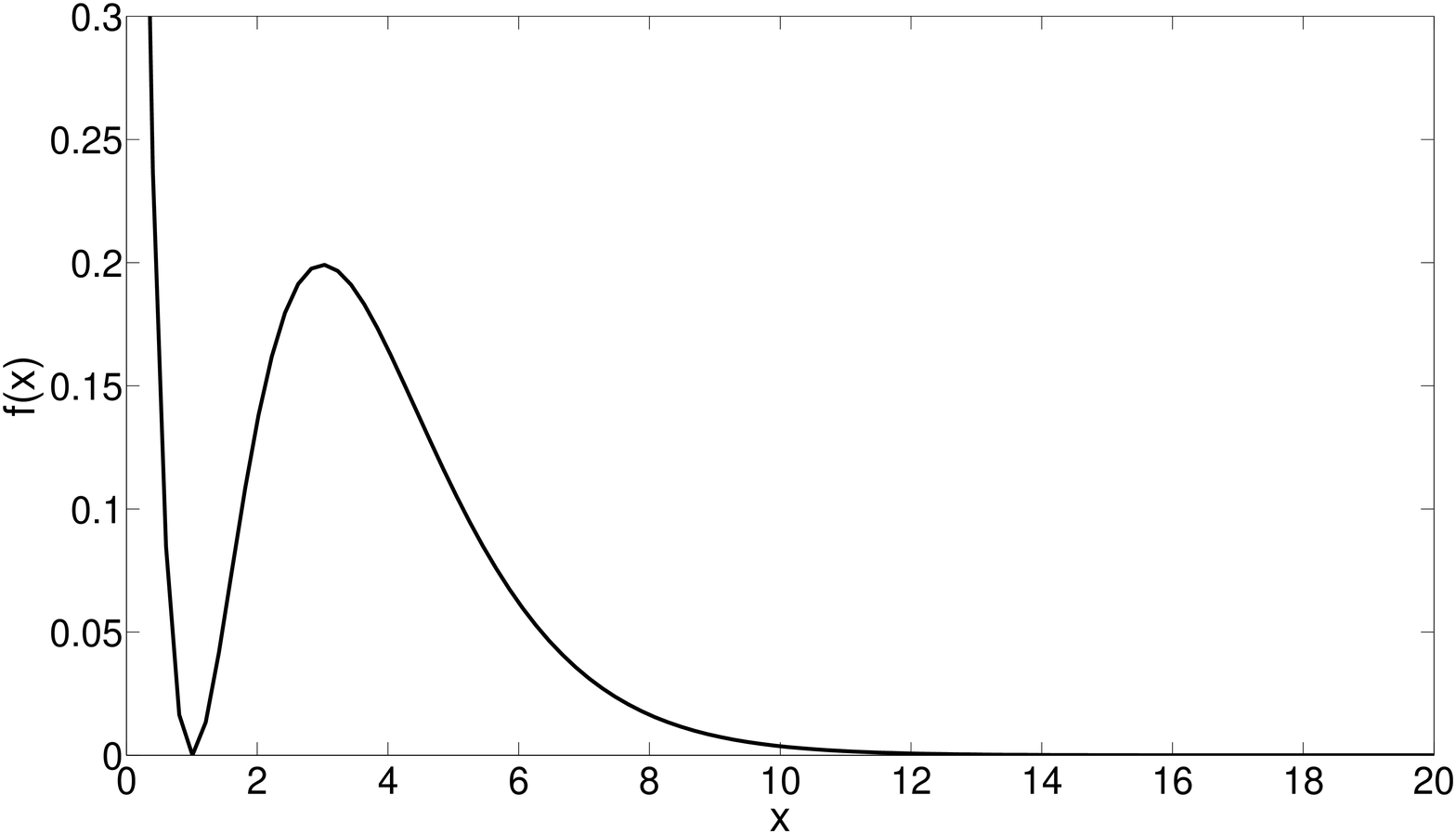}\hspace{0.5cm}
\includegraphics*[width=0.5\textwidth]{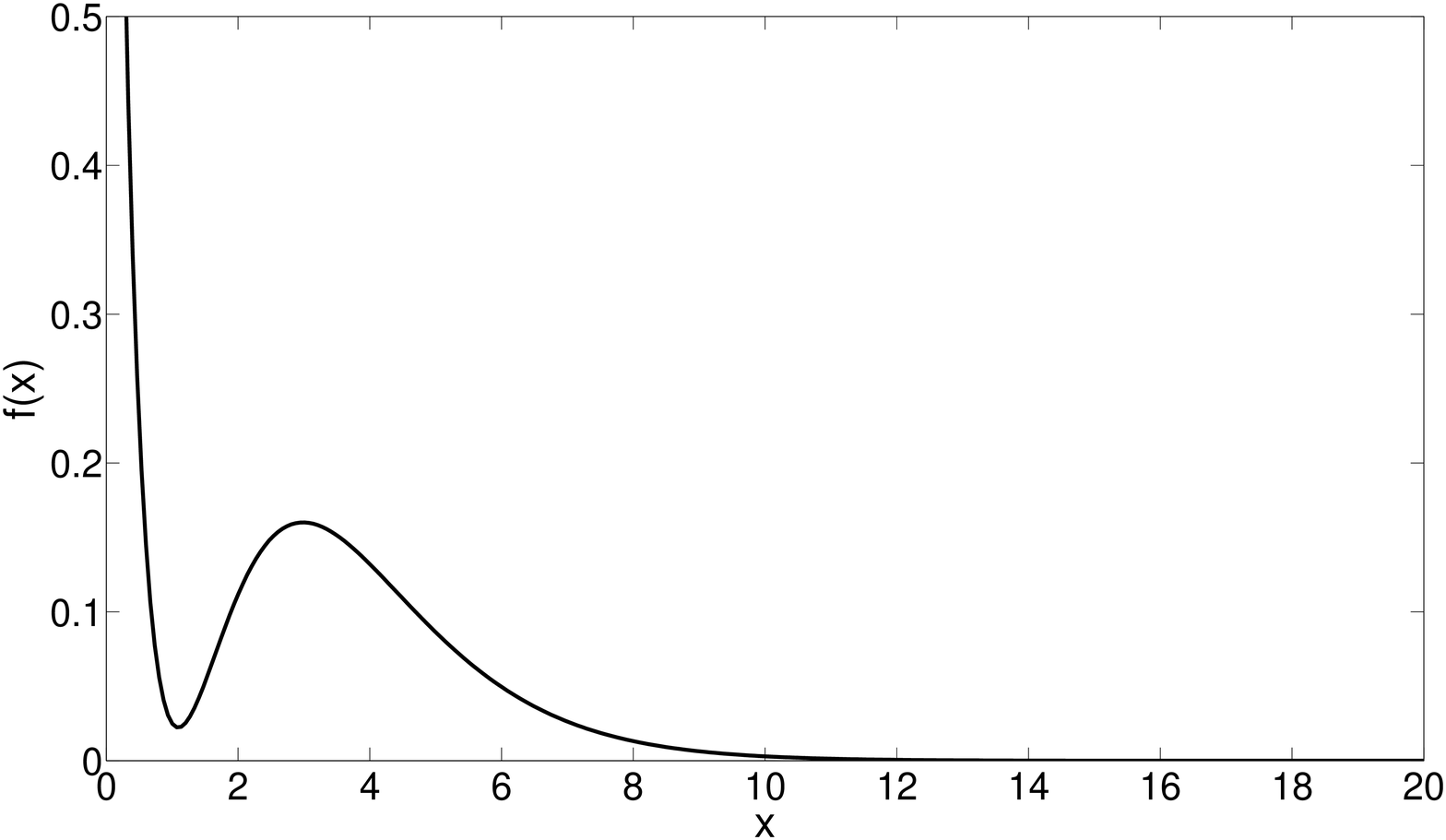}
\caption{Initial distributions in the example d) (left) and  the example e) with $\alpha=2$ (right)}
\label{F0}
\end{figure}

{\bf d)} Finally let us show that the coarsening rates  for large times depend
only on how fast the initial distribution $f(x)$ decays as $x\go+\infty$
and not its behavior for moderate $x$. In this and the next example we consider non-monotone distributions having a local maximum at $x>0$.
Consider $f(x)=(1-x)^2\exp(-x)$ (see Fig. \ref{F0}) with $n(x)=(1+x^2)\exp(-x)$, correspondingly. By \rf{CRL} one obtains
\bea
\hspace{-0.7cm}BT(d)&=&\int_0^d(1+x^2)\exp(-x)\ln(1+x^2)\,dx-\int_0^dx(1+x^2)\exp(-x)\,dx\nonumber\\ &-&(\ln(1+d^2)-d)\int_0^d(1+x^2)\exp(-x)\,dx\nonumber\\
&=&\int_0^d(1+x^2)\exp(-x)\ln(1+x^2)\,dx-7+3(d-\ln(1+d^2))+O(\exp(-d)).
\lb{Ms8}
\eea
The first integral in the last expression can be estimated as follows.
\beas
\int_0^d(1+x^2)\exp(-x)\ln(1+x^2)\,dx&\le&\ln(1+d^2)\int_0^d(1+x^2)\exp(-x)\,dx\\
&=&\ln(1+d^2)(3+O(\exp-d))
\eeas
Combining this with \rf{Ms8} one obtains
\bes
T(d)=\frac{3d}{B}+o(d)
\ees
and hence the following asymptotics holds
\bes
n(t)\sim \left(1+\frac{9}{B^2}t^2\right)\exp(-Bt)\ \ \text{as}\ \ t\go\infty.
\ees
Therefore, the coarsening rates show an exponential decay as in the example b).

{\bf e)} Consider distributions
\be
\frac{\alpha}{\alpha+1}\left[(1-x)^2\exp(-x)+1/(1+x)^{1+\alpha}\right]\ \ \text{with}\ \ \alpha>0,\,\alpha\ne 1.
\lb{PDF1}
\ee
They have a local maximum at $x>0$ and a decay $\sim 1/x^{1+\alpha}$ as $x\go \infty$ (see  Fig. \ref{F0}).
Correspondingly, one has
\bes
n(x)=\frac{1}{1+\alpha}\left[(1+x)^{-\alpha}+\alpha\exp-x(1+x^2)\right].
\ees
Substituting it in \rf{CRL} one obtains
\beas
BT(d)&=&\int_0^d\left[\int_0^xn(y)\,dy\right]\frac{n'(x)}{n(x)}\,dx=\\
&=&\int_0^d\left[\frac{1}{1-\alpha}\left((1+x)^{1-\alpha}-1\right)+\alpha\left(3-\exp(-d)(3+d(2+d))\right)\right]\times\\
&\times&\frac{(1-x)^2\exp(-x)+(1+x)^{-1-\alpha}}{(1+x)^{-\alpha}+\alpha\exp(-x)(1+x^2)}\,dx.
\eeas
The last integral can be bounded from below and from above by integrals of the following type
\bes
I^*=\int_0^d\left[\frac{1}{1-\alpha}\left((1+x)^{1-\alpha}-1\right)+C_1\right]
\frac{C_2\exp(-x/2)+(1+x)^{-1}}{1+C_3}
\lb{Ms9}
\ees
with some nonnegative constants $C_i,\ i=1,2,3$. Integrals in \rf{Ms9} have the following asymptotics:
\bes
I^*=C_3\left(\frac{1}{1-\alpha}\left[(d+1)^{1-\alpha}-1\right]-C_4\ln(d+1)\right)+O(1)\ \ \text{as}\ \ d\go\infty
\ees
with some positive constants $C_3,\,C_4$. 
Therefore, using the asymptotics
\be
n(d)\sim \frac{1}{1+\alpha}(1+d)^{-\alpha}\ \ \text{as}\ \ d\go\infty
\lb{Ms10}
\ee
one obtains that the asymptotics of the coarsening law for \rf{PDF1} 
coincides up to multiplicative constants with \rf{alpha} already obtained in the example a) for the monotone distributions.
Consequently, we conclude that as in the example a) the coarsening rates are algebraic with power $\alpha/(\alpha-1)$ for $\alpha<1$ and 
exponential for $\alpha>1$. 

A more simple but rather formal proof of this fact is as follows. Let us fix a large number $A$ such that asymptotics \rf{Ms10} holds for all $d>A$ with
a good precision. The one has
\beas
BT(d)&=&\int_0^An(x)\ln\left[\frac{n(x)}{n(d)}\right]\,dx+\int_A^dn(x)\ln\left[\frac{n(x)}{n(d)}\right]\,dx\\[2ex]
&\sim&O(1)+\alpha\ln(d+1)\times O(1)+\int_A^dn(x)\ln\left[\frac{n(x)}{n(d)}\right]\,dx.
\eeas 
The last integral in view of \rf{Ms10} is of the type considered already in
example a). Therefore,  the term $O(1)+\alpha\ln(d+1)\times O(1)$ produces no change in the asymptotics of
 $T(d)$ and, consequently, the coarsening law coincides up to multiplicative constants with \rf{alpha}.

\section{Numerics} 

In the last two sections starting from the system \rf{SSM1}-\rf{SSM2} we first derived a closed ODE model \rf{ODE}--\rf{PB}
describing coarsening dynamics in an array of initially $N+1$ metastable droplets. Then we looked at its limiting case $\beta\go \infty$ 
described by \rf{ODEinf}--\rf{PBinf} and more precisely on its leading order
version as $\eps\go 0$ given by \rf{RS}. Next, we found out that for a special initial
data satisfying \rf{ID} one can obtain the explicit solution to \rf{RS} given by \rf{RSS}. Assuming additionally \rf{ID2} and that the distances
in the array are ordered decreasingly we derived an explicit coarsening law \rf{DL}. Finally, we obtained its continuous counterpart \rf{CRL}.
In this section we systematically compare numerically the solutions of subsequent models in the derived model hierarchy and check coarsening laws 
\rf{DL}, \rf{CRL}.

\subsection{Comparison between models} 
Here we compare solutions to the full ODE system \rf{ODE}--\rf{PB} with those of the strong-slip system \rf{SSM1}-\rf{SSM2} and its limiting cases
\rf{ISM} and \rf{FFM1}-\rf{FFM2} as $\beta\go 0$ and  $\beta\go +\infty$, respectively. For the solution of PDE systems we used a fully implicit finite difference scheme derived and applied already to \rf{SSM1}-\rf{SSM2} and its limiting cases in~\cite{muench04a, MWW06,peschka08,Ki09}. 
The numerical solutions for \rf{ODE}--\rf{PB} were obtained applying a fourth-order adaptive time step 
Runge-Kutta method and using updating rules \rf{Pmax}--\rf{ClnRecal} after each subsequent coarsening event. In the case of PDE system the corresponding pressure evolution was calculated using finite-difference discretization of the term $\Pi_\eps(h)-\pxx h$.
\begin{figure}
\includegraphics*[width=0.4\textwidth]{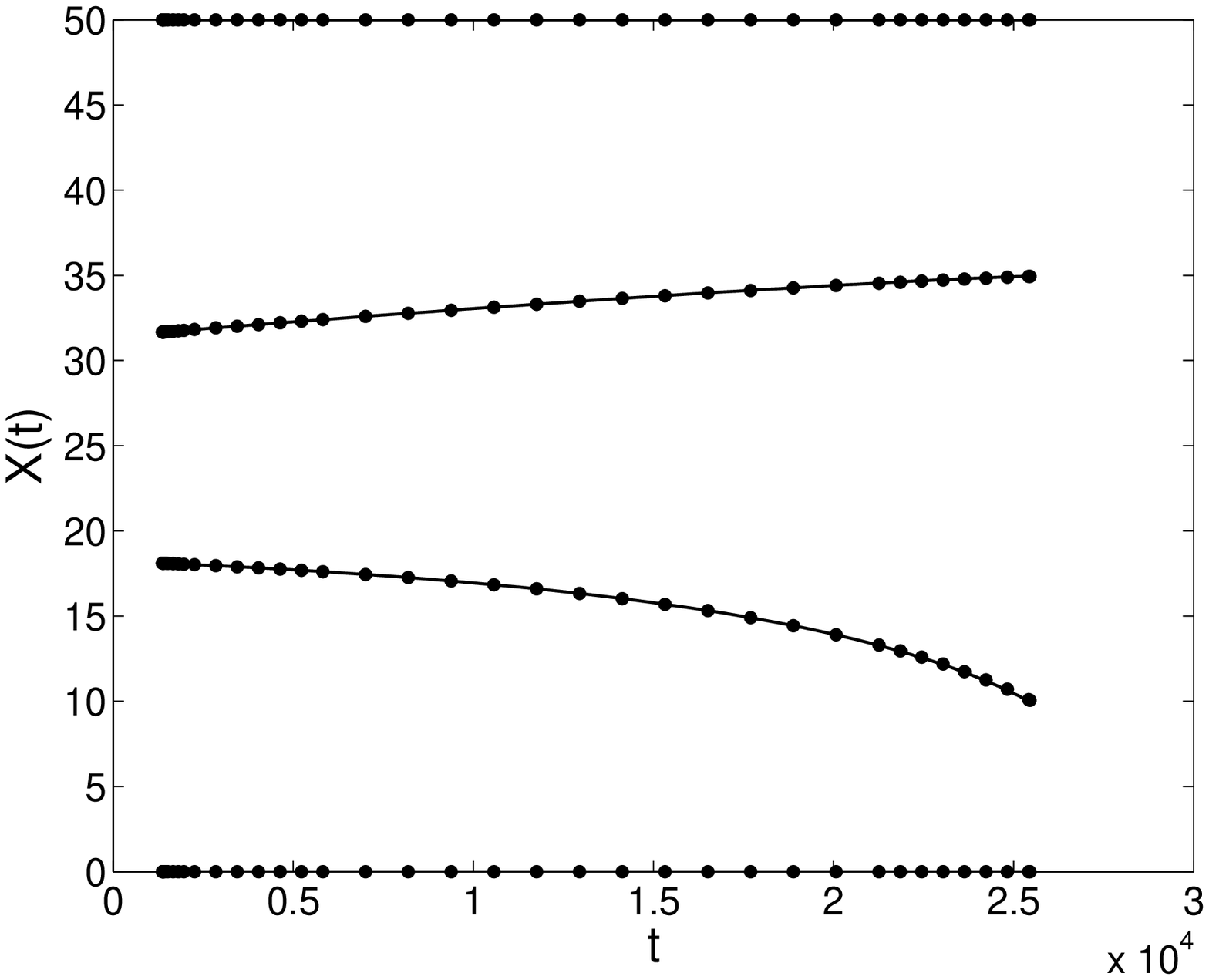}\hspace{2.5cm}
\includegraphics*[width=0.4\textwidth]{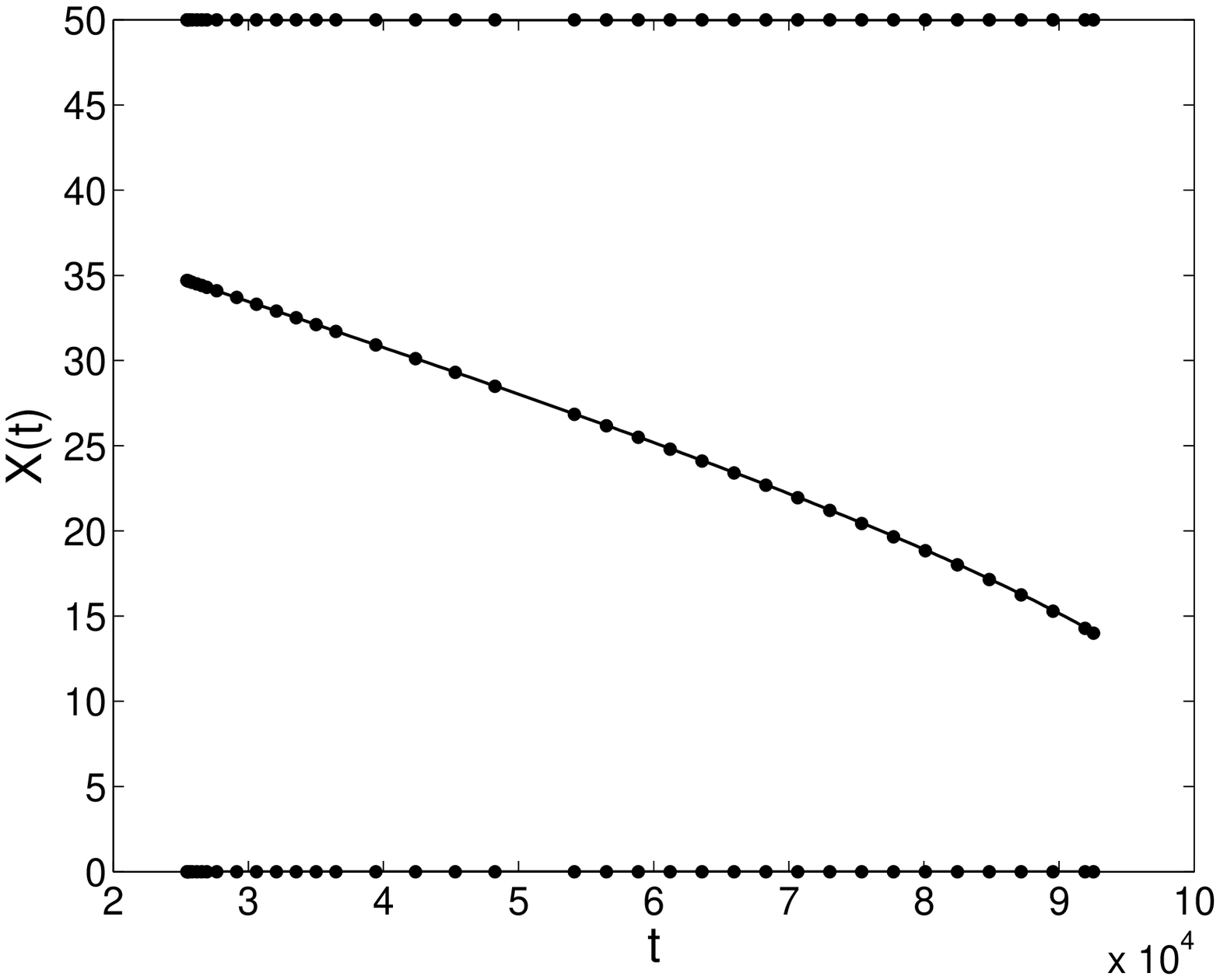}
\caption{Comparison of droplet position evolution obtained from the ODE model (dots) and the 
lubrication model (solid line) in the strong-slip case with $\eps=0.025$,
$Re=0$, $\beta=10$. Initial profile of four droplets presented on the top-left
plot of Fig. \ref{F2} was used. Subsequent collisions of the second (left) and the third (right) droplets with
the first one are shown.}
\label{F1}
\end{figure}

In Fig. \ref{F1} starting from an array of four droplets we compare evolution of
positions resulting from PDE and ODE models for two subsequent collisions. 
Fig. \ref{F1} shows that the absolute deviation between results stays uniformly $O(\eps)$
also after subsequent collision events. 

\begin{figure}
\includegraphics*[width=0.4\textwidth]{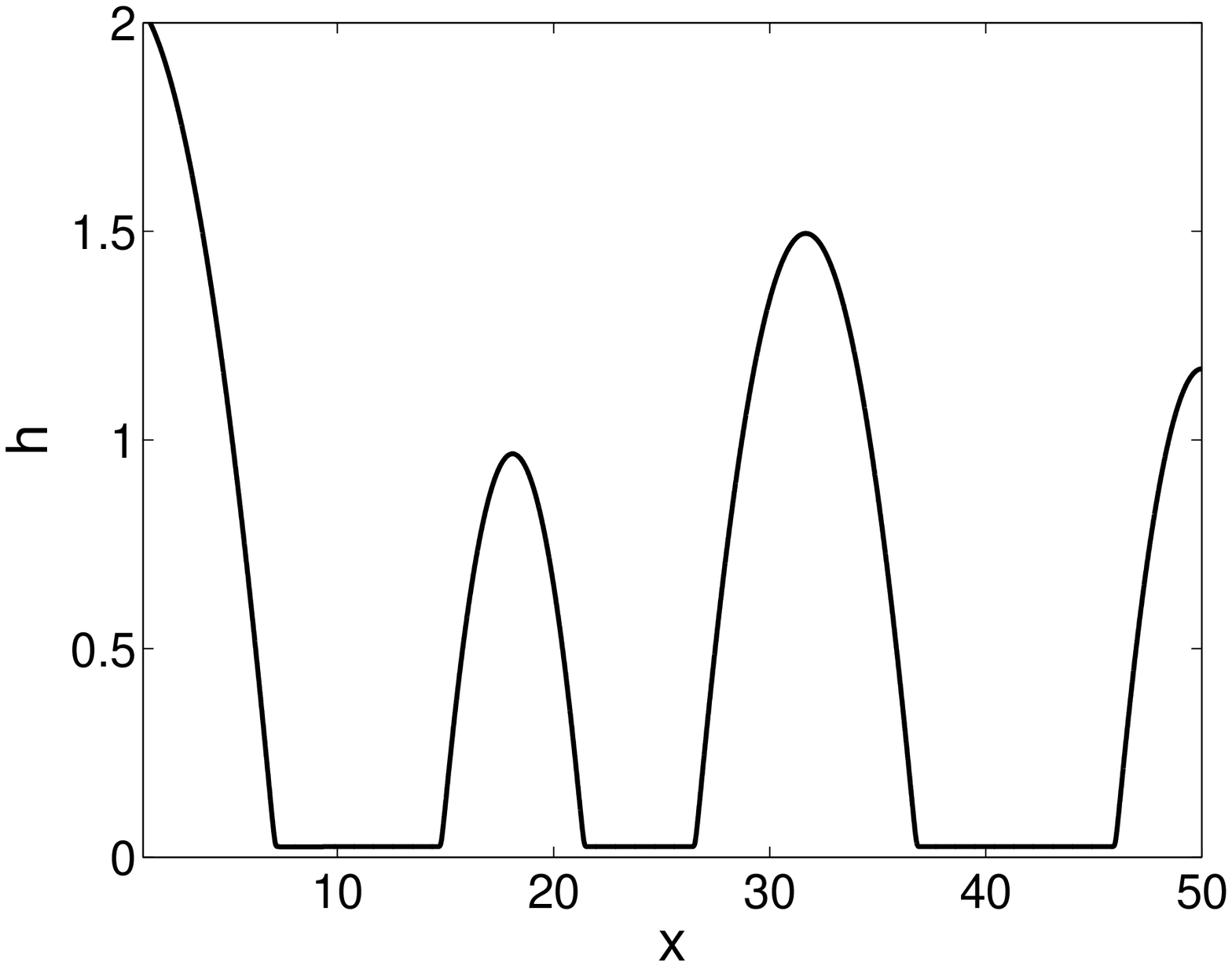}\hspace{2.5cm}
\includegraphics*[width=0.4\textwidth]{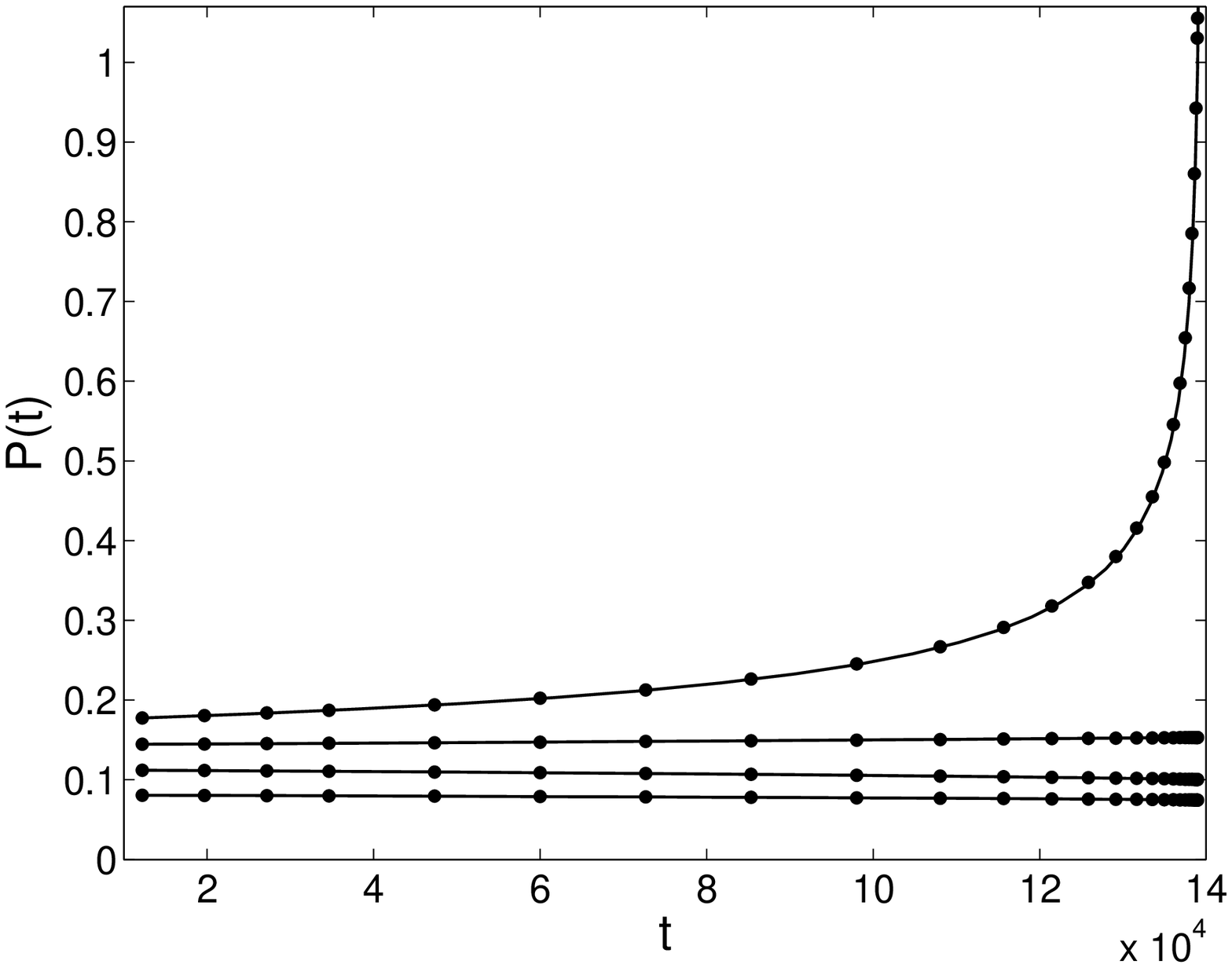}

\includegraphics*[width=0.4\textwidth]{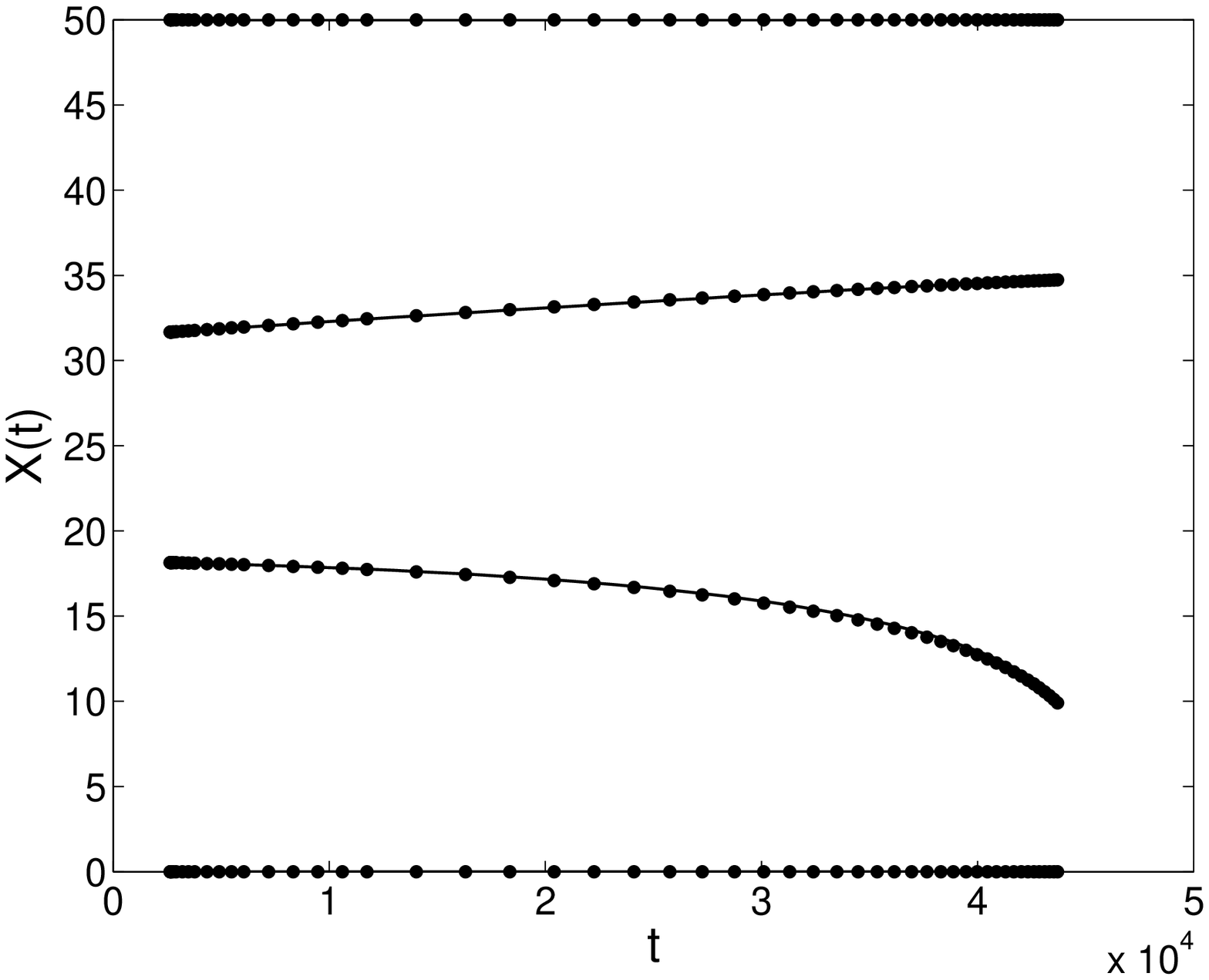}\hspace{2.5cm}
\includegraphics*[width=0.4\textwidth]{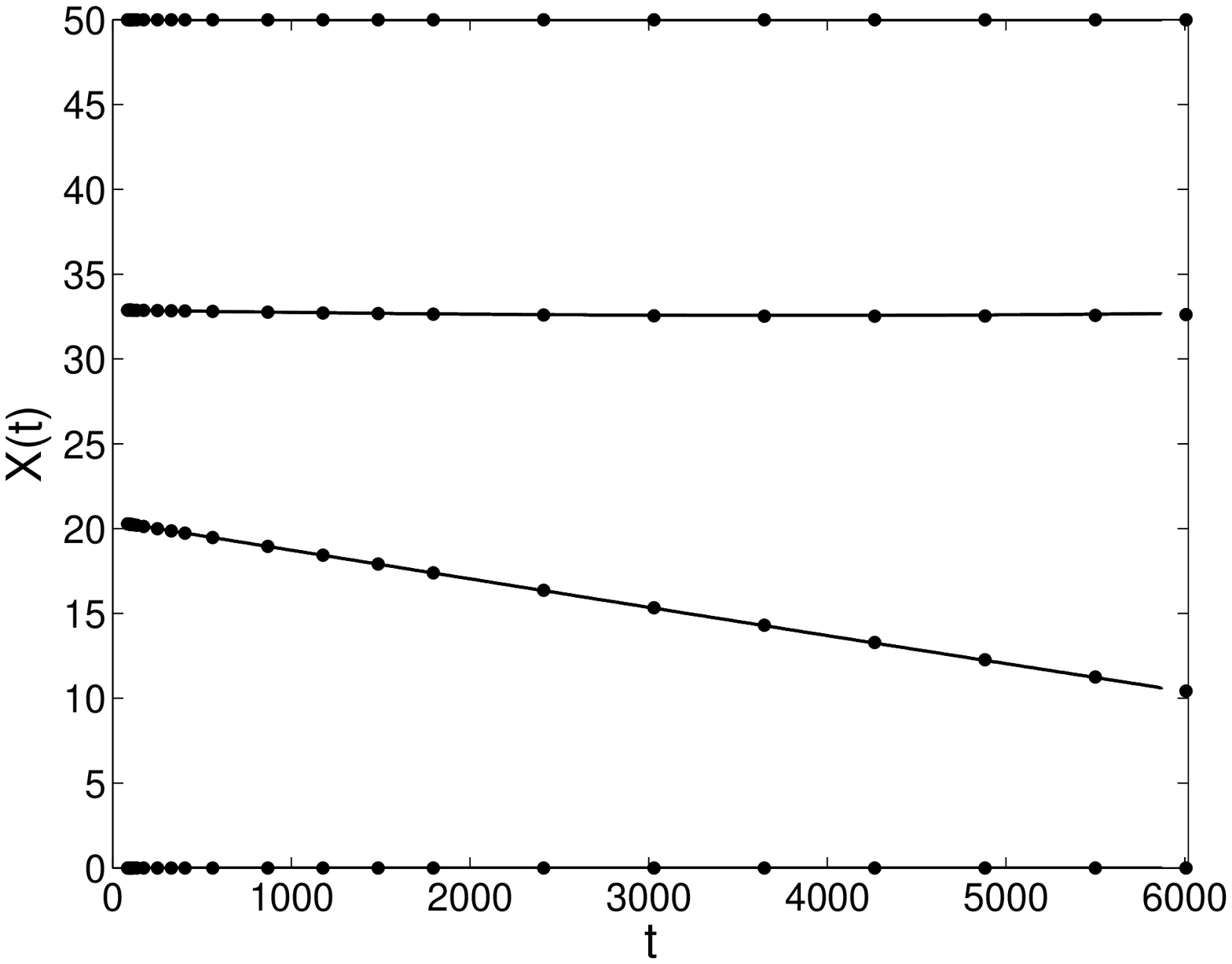}
\caption{Comparison of droplet evolution obtained from the ODE model (dots) and the 
lubrication model (solid line) in the strong-slip case with $\eps=0.025$,
$Re=0$, and different $\beta$. Upper row: Initial profile of four droplets
(left) and their pressure evolution in the intermediate-slip case until collapse of the second droplet (right).
Lower row: Starting from the same initial profile position evolution until collision of the second droplet with the
first one for $\beta=5$ (left) and $\beta=\infty$ (right) is shown.}
\label{F2}
\end{figure}

In Fig. \ref{F2}  starting from the same array of four droplets we compare
solutions for different slip lengths $\beta$. 
In the cases $\beta=0$ and $\beta=\infty$ we compared solutions to \rf{ISM}
and \rf{FFM1}-\rf{FFM2} and those of \rf{ODE0}--\rf{PB0} and
\rf{ODEinf}--\rf{PBinf}, respectively. Again for all $\beta$ the absolute
deviation between PDE and ODE results is $O(\eps)$.
Note, Fig. \ref{F2} demonstrates the fact pointed already in~\cite{ORS07,KW09} that in
the intermediate-slip case the coarsening dynamics is governed mostly by
collapse mechanism while in the strong-slip case with moderate and large
$\beta$ by collisions.

\begin{figure}
\includegraphics*[width=0.4\textwidth]{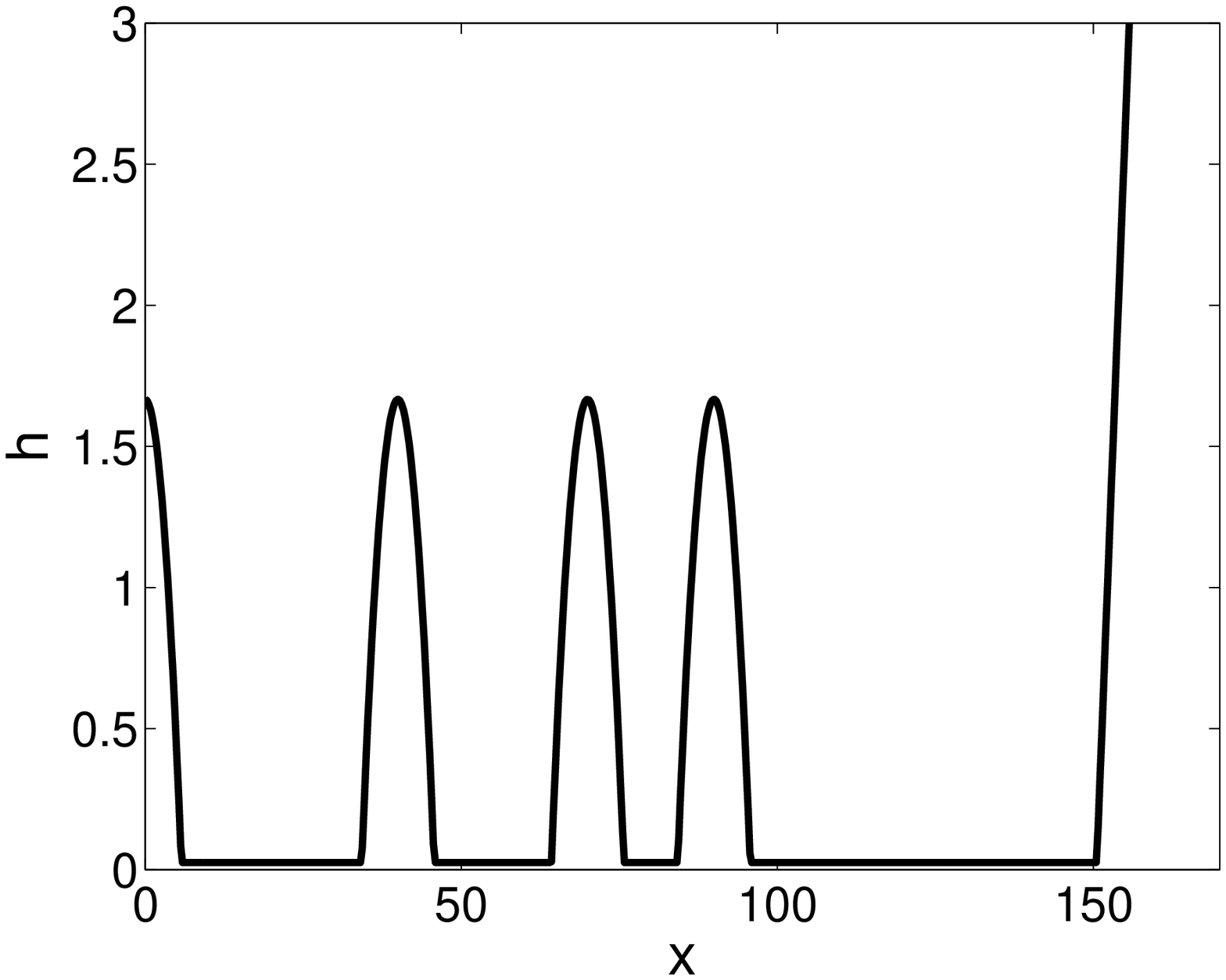}\hspace{2.5cm}
\includegraphics*[width=0.4\textwidth]{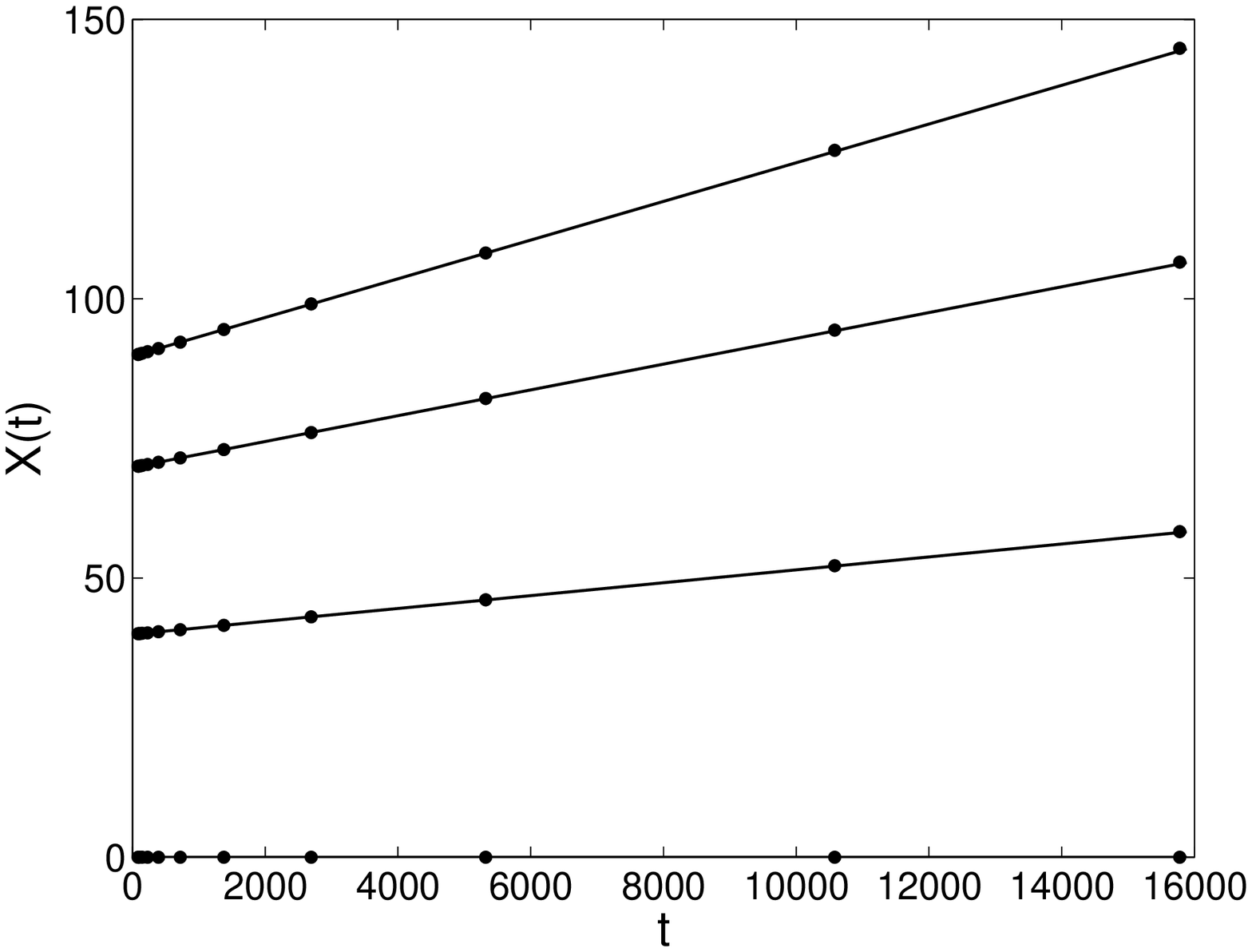}
\caption{Comparison of droplet position evolution obtained from the ODE model
\rf{ODEinf}--\rf{PBinf} (dots)with $\eps=0.025$ and its leading order
subsystem \rf{RS} (solid line).  
Left plot: The initial profile of five droplets. Pressure of the first four
and the last droplet are $0.01$ and $0.001$, respectively. Right plot: Collision of the fourth droplet with the largest last one.}
\label{F3}
\end{figure}

Finally,  Fig. \ref{F3} shows that for arrays being taken initially with sufficiently
small pressures the migration subsystem \rf{RS} approximates with a high
accuracy (at least $O(\eps)$) the full ODE system \rf{ODEinf}--\rf{PBinf}
describing  the case $\beta=\infty$. Note, that the migration path of the fourth droplet Fig. \ref{F3} is
considerably large $\approx 50$. Hence, the reduction from
\rf{ODEinf}--\rf{PBinf} to \rf{RS} does not constrains the initial distances to be small.

\subsection{Coarsening rates} 
\begin{figure}
\includegraphics*[width=0.4\textwidth]{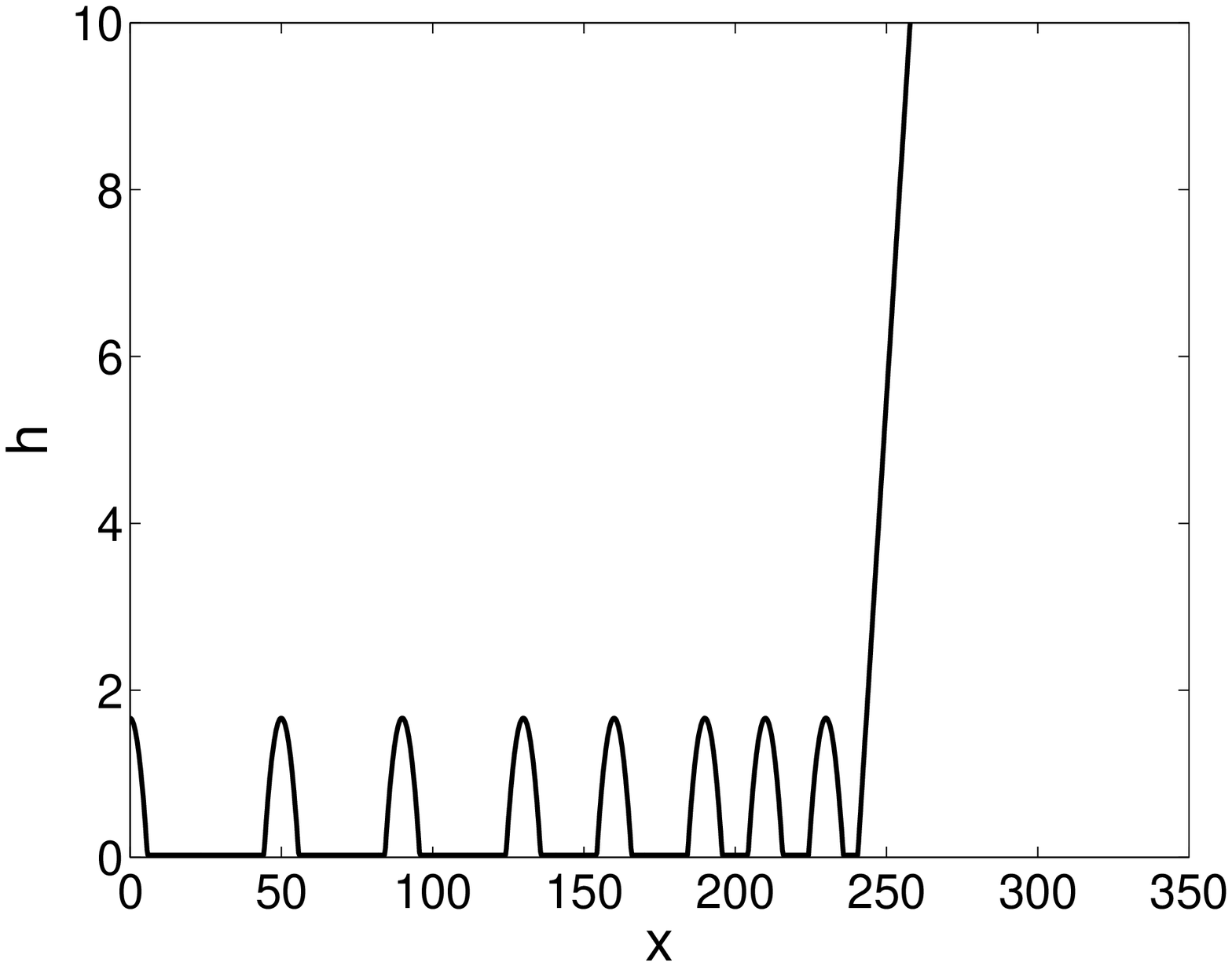}\hspace{2.5cm}
\includegraphics*[width=0.4\textwidth]{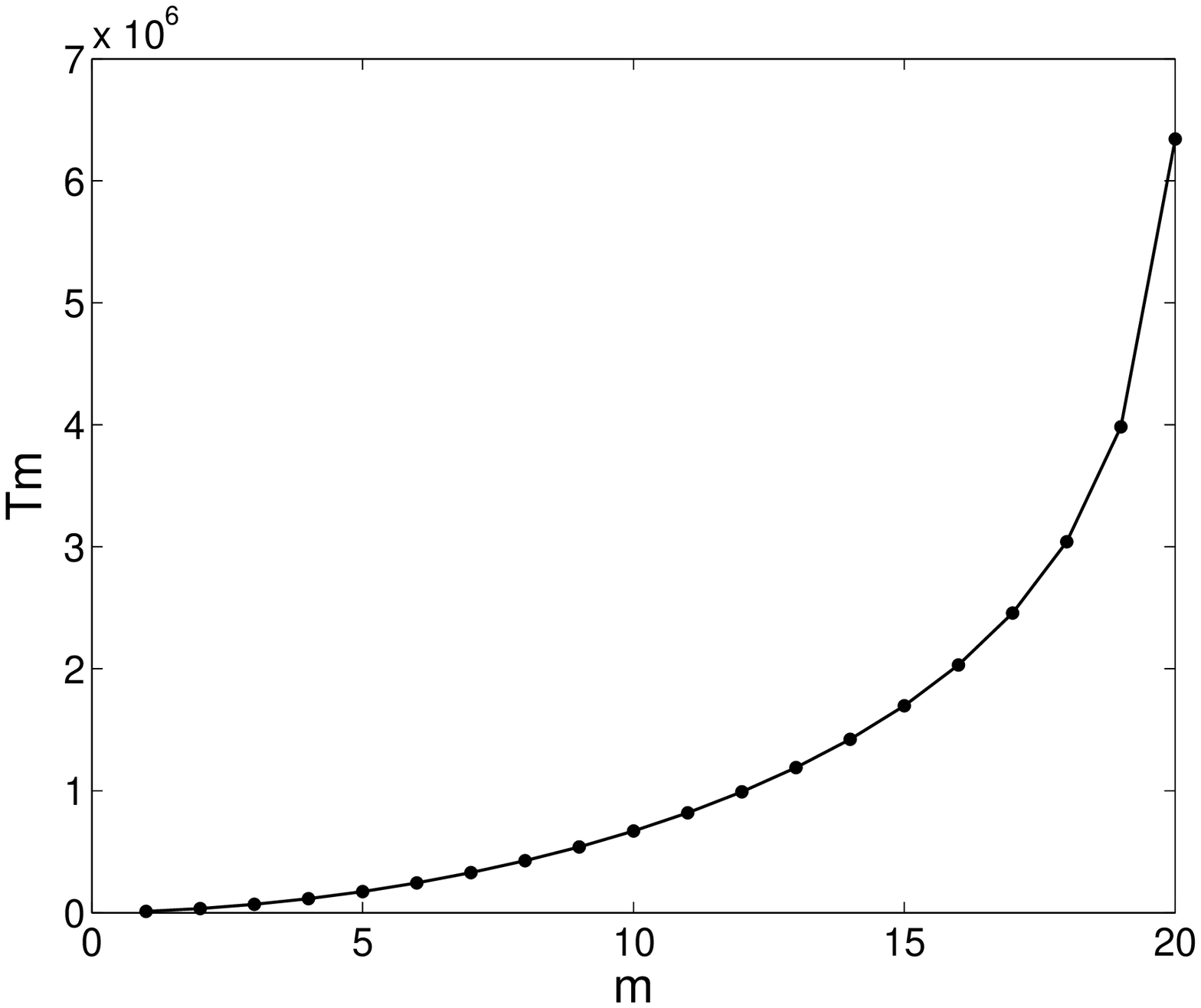}
\caption{Comparison of the subsequent
collision/absorption times for each $m$-th family ($1\le m\le 20$) given by the discrete
coarsening law \rf{DL} (dots) and iterative calculation using the solution \rf{RSS} (solid line).
Left plot: A part of a typical initial profile under consideration. Right plot: plot of
the absorption times versus the family number.}
\label{F4}
\end{figure}

Here using the explicit solution \rf{RSS} for the system \rf{RS} we
check numerically the discrete and continuous coarsening laws \rf{DL} and \rf{CRL}.
In Fig. \ref{F4} we take $k=20$ families of initial distances as prescribed in
\rf{hi} with the corresponding pressures satisfying \rf{ID}, \rf{ID2}
and order them non increasingly in space. Next, we compare the subsequent
times for each $m$-th family ($1\le m\le k$) to be absorbed given on one hand by the 
analytical law \rf{DL} and by iterative calculation using \rf{RSS},\rf{hu}
between each subsequent collision  on the other. Naturally, one finds out the exact coincidence between them.
Note, as each collision is comprised of an absorption by the largest droplet of a smaller one one do not need
to update the position and pressure of the former one. This is because its position is fixed due to \rf{BC} to $x=L$ and the pressure to the leading order does not change due to \rf{ID2}.

\begin{figure}
\includegraphics*[width=0.4\textwidth]{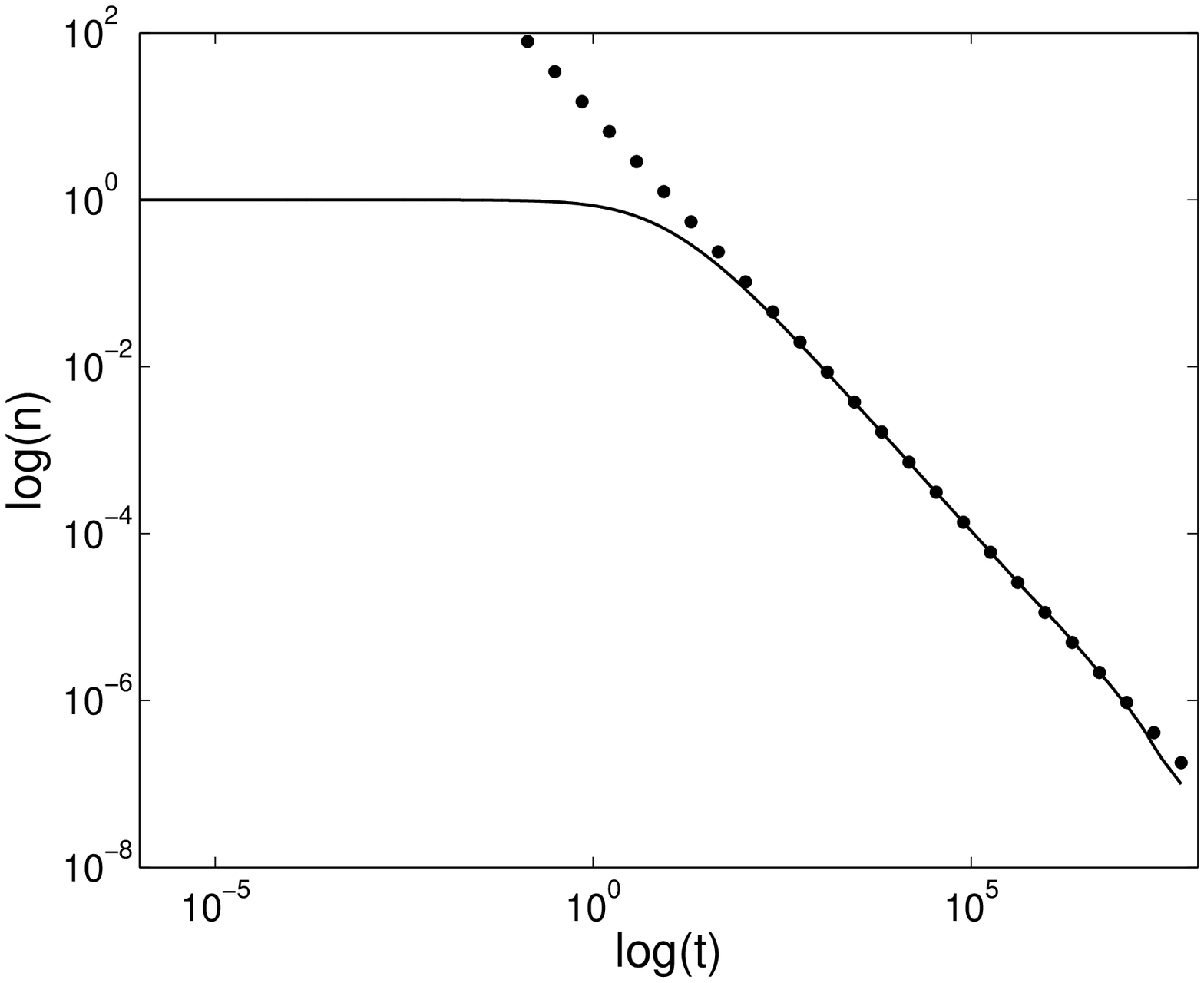}\hspace{2.5cm}
\includegraphics*[width=0.4\textwidth]{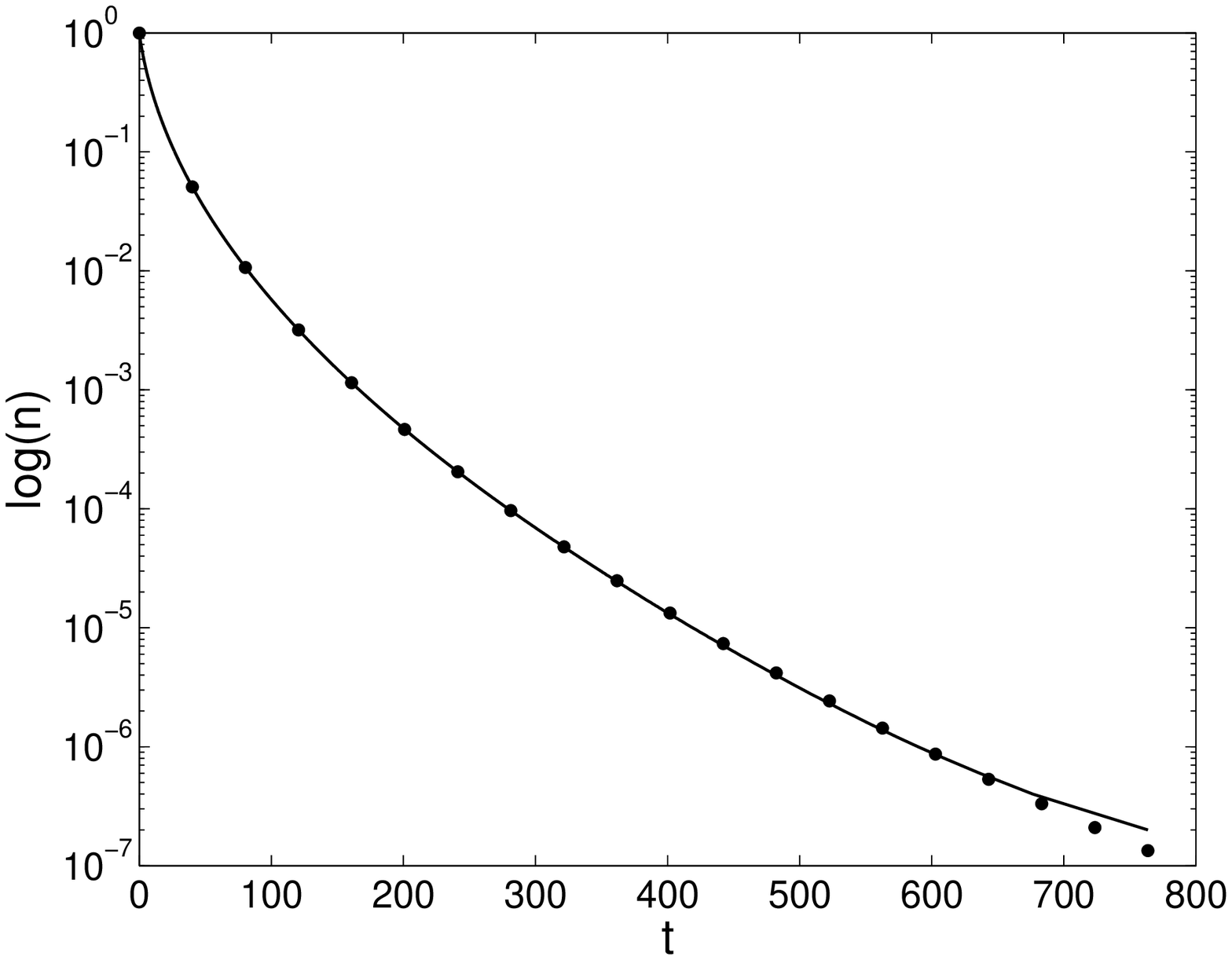}

\includegraphics*[width=0.4\textwidth]{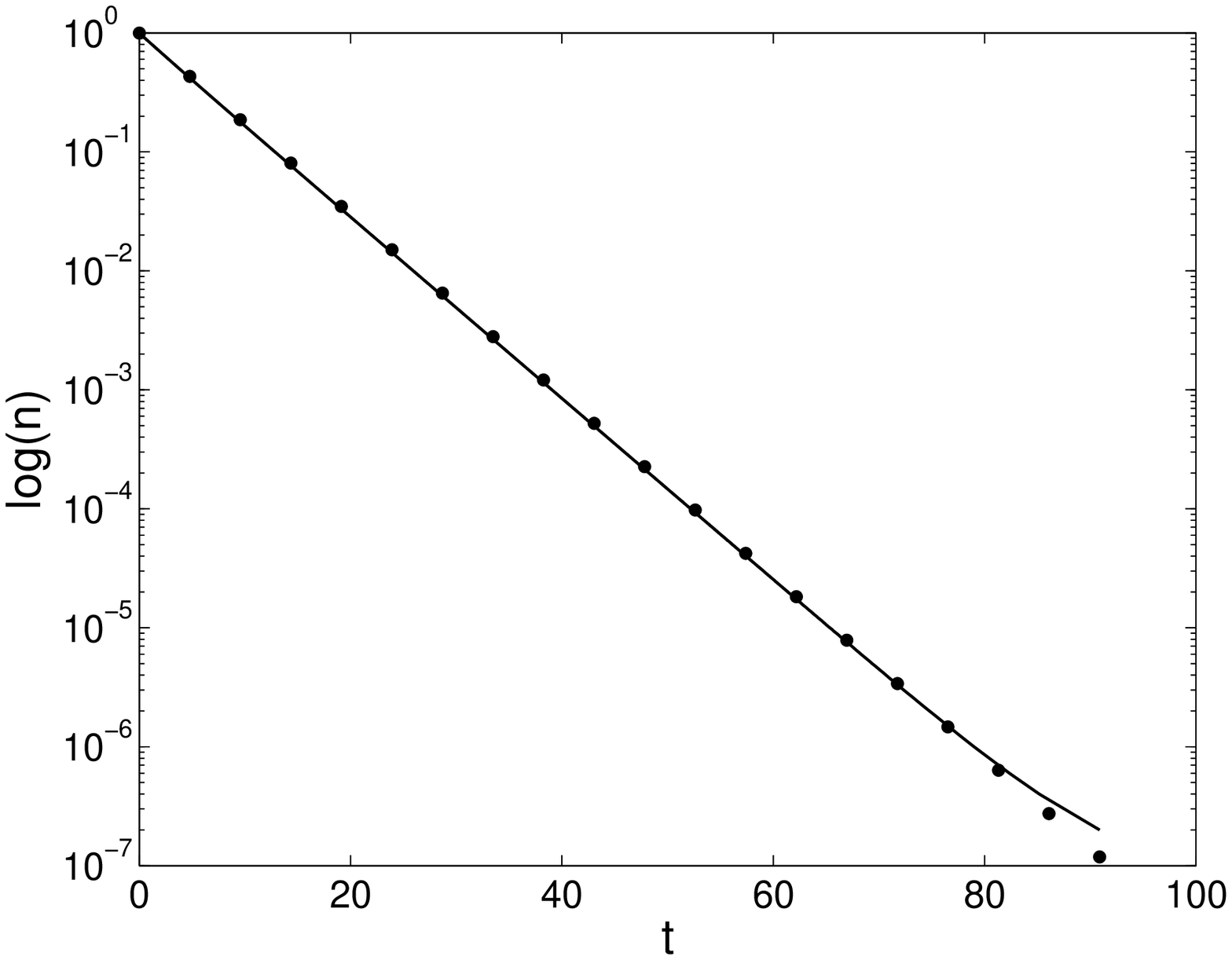}\hspace{2.5cm}
\includegraphics*[width=0.4\textwidth]{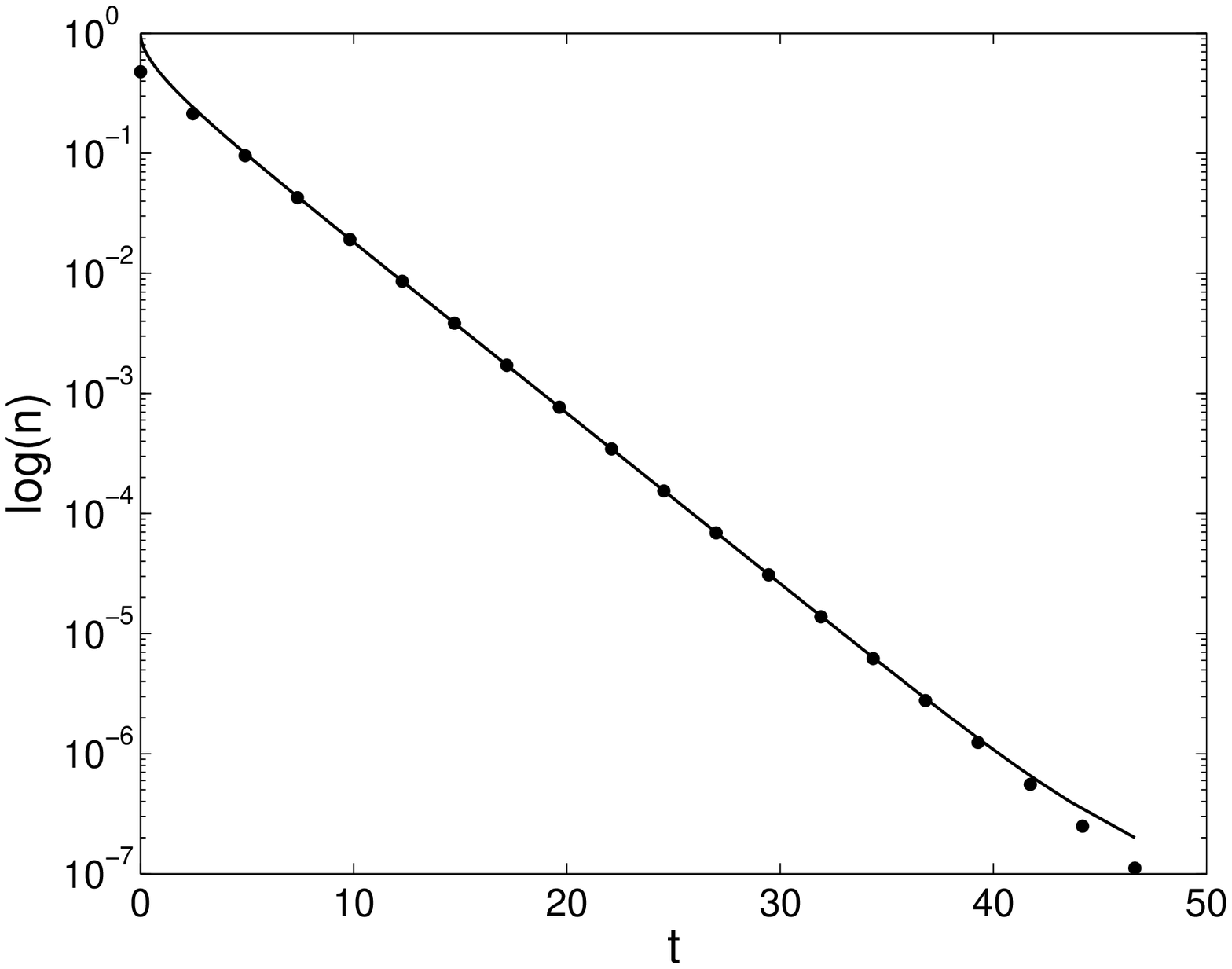}
\caption{Comparison of numerical coarsening rates using sampling of
the initial data with $N\gg 1$ and subsequent iterative calculation using \rf{RSS} (solid line) with those provided by the analytical law \rf{CRL} (dots).
Upper row: loglog and semilog plots for the initial distribution \rf{PDF} with $\alpha=1/2$
(left) and $\alpha=1$ (right), respectively.
Lower row: semilog plots for \rf{PDF} with $\alpha=20$
(left) and  for the Gaussian initial distribution (right). According to \rf{CR} and \rf{CRgauss} in the chosen axe scales
the dots reproduce the linear functions except for the upper-right plot where they represent the law \rf{CR1}}
\label{F5}
\end{figure}

In Fig. \ref{F5} we show numerical results for continuous coarsening rates
for three initial distributions taken from the family \rf{PDF} with different $\alpha$ and one
Gaussian distribution considered in examples a) and c) of the previous section, respectively. 
To obtain the coarsening rates numerically we first sampled $N\gg 1$ distances
according to the given initial distribution. After ordering them non increasingly
in the initial configuration we substitute them as initial data into \rf{RS} and solve the latter one iteratively
using \rf{RSS}.  Note, that due to an extremal simplicity of \rf{RSS}
one can effectively model numerically a huge number of droplets $N\approx 10^7$ just using capabilities
of a personal computer. 

Fig. \ref{F5} shows that thus obtained numerical coarsening rates coincide for large
times very well with the analytical ones prescribed by the law \rf{CRL} and found out in examples a) and c) of the section 4.
In the case of \rf{PDF} with $\alpha=1$ one has the exact coarsening law
\rf{CR1} and therefore a good coincidence for all times. In the case of
\rf{PDF} with $\alpha\ne 1$ we compared our numerical results with the asymptotic law \rf{CR},
while for the Gaussian initial distribution with \rf{CRgauss}. Note that a certain deviation between numerics and analytical laws starting 
at the very end of the considered time interval is caused by a numerical error increase in sampling of large distances according to 
initial probability distributions.

\section{Conclusions and discussion} 
In this paper we started from the lubrication equations \rf{SSM1}--\rf{SSM2} describing dewetting process in nanometric polymer film
interacting on a hydrophobically coated solid substrate in the presence of large slippage at the liquid/solid interface.
This model describes a distinguished and important regime within a lubrication scaling. 
In particular, it incorporates as a limiting case of the infinite slip length 
the  well-known model of free films \rf{FFM1}--\rf{FFM2} studied intensively in applications~\cite{ED93,KSK04}.
Note, that a similar to \rf{FFM1}--\rf{FFM2} system appears in the study of viscoelastic threads 
for which coarsening dynamics of interacting droplets was observed also at the experiments, see e.g. \cite{CELM06}. 

Motivated by this we derived the reduced ODE models \rf{ODE}-\rf{PB0}
describing coarsening dynamics of droplets governed by \rf{SSM1}--\rf{SSM2} and its limiting cases.  In the limiting case
$\beta=\infty$ we observed that the migration subsystem \rf{RS} can be solved explicitly for special initial data satisfying \rf{ID}, \rf{ID2}.
By \rf{RSS}-\rf{hu} the dynamics of droplets consists of sequential collisions of smaller ones with the largest one while their distance is distributed
uniformly between the remaining drops. 

Similar models were suggested for collapse/collision dynamics of breath figures by~\citet{DGY91}
basing on heuristic arguments. There authors considered the 'cut-in-two' and 'past-all' models where the distance of the smallest droplet was divided between two neighbors or pasted as a whole to one of them, respectively.  
The breath figures of~\cite{DGY91} found later interesting analogs in the reduced coarsening models 
arising from  the Allen-Cahn and Ginsburg-Landau equations,see e.g. a review article~\cite{B94}. A recent generalization of these models and rigorous analysis
of their self-similar solutions can be found in~\cite{MBP10}. Note that \rf{RS} can be classified due to \rf{hu} as a 'cut uniformly between all' type model.
We are not aware if any heuristic or rigorous analog of it was considered so far in the literature.

Remarkably, our model \rf{RS} appears as a reduction of a complicated dynamics governed 
by a high-order lubrication system \rf{FFM1}--\rf{FFM2}. Moreover, in contrast to stochastic in nature models of~\cite{DGY91}
if the initial distances in \rf{RS} are ordered non increasingly  then they coarse in a deterministic and exactly solvable scenario.  
For the latter case we derived the coarsening laws \rf{DL} and \rf{CRL} analytically and confirmed  them numerically. 
Interestingly, the derived law  \rf{CRL} have a form similar to the  Shannon's entropy with respect to a  certain
normalization of the droplet number function $n(x)$. The explanation of this fact from the statistical point of view will be presented elsewhere.

Surprisingly, in contrast to the coarsening dynamics governed by reduced ODE models arising from \rf{WSM} which obey always the law \rf{NCL}
our simple model \rf{RS} can reproduce any algebraic coarsening rates between zero and infinity as well as exponential ones.
Moreover, for a family of initial distributions \rf{PDF} we showed existence of a threshold for their decay at infinity 
at which the corresponding coarsening rates switch from algebraic to exponential ones. 
Note, that a similar situation was accounted recently for self-similar solutions to
Smoluchowski coagulation equation with certain kernels, see~\cite{MP04,MP08}. 

In view of the above observations it would be natural to extend our deterministic collision/absorption model
to its stochastic variant withdrawing  the non increasing order of the distances and thus allowing collisions of random droplets with the largest one.
As in~\cite{DGY91} one could  probably look for self-similar solutions of the mean-field approximations for thus arising stochastic collision models.
A further generalization of the model could be a withdrawal of the constraints on the initial data \rf{ID}, \rf{ID2}  and thus allowing droplets to collide and collapse inside of the domain. Note, that an additional difficulty to handle the pressure and position update according to the coarsening rules
\rf{Pmax}--\rf{ClnRecal} would appear then.

Finally, it could be possible to derive two-dimensional analogs
of the reduced ODE models \rf{ODE}-\rf{PB0} describing physically coarsening of three dimensional droplets on a plane substrate.
The two dimensional reduced ODE models arising from \rf{WSM} were derived in~\cite{Gl07,ORS07}. In~\cite{Gl07}
a mean-filed approximation for the fluxes between droplets was suggested under an assumption of well separation of droplets 
that is unfortunately not suitable for the modelling of droplet collisions because the distance between them tends to zero then. 
In this case one should face a problem of solving a Laplace equation (counterpart to equation \rf{j2} in two dimensions)
in a complex domain between droplets occupied by the PL region. We expect the same problem to appear for the two-dimensional
reduced ODE models corresponding to \rf{SSM1}--\rf{SSM2}.   

\section{Acknowledgments} 

The author acknowledges the postdoctoral scholarship at the Max-Planck-Institute for Mathematics in the Natural Sciences, Leipzig
and thanks Gennady Chuev, Christian Seis and Andre Schlichting for fruitful discussion.

\begin{appendix}
\section{Integral I} 
Here we show that integral $I$ defined in \rf{I} converges and integrate it explicitly.
Changing variable in \rf{I} according to the explicit solution \rf{H1} in the CL region, 
and using matching conditions \rf{MC1} one obtains
\beas
I&=&\int_{-\infty}^{+\infty}\frac{1}{H_1}\pz\left(\frac{\pz H_1}{H_1}\right)\,dz=\int_1^{+\infty}\left[\frac{U'(H)}{H^2}-\frac{2(U(H)-U(1))}{H^3}\right]
\frac{dH}{\sqrt{2(U(H)-U(1))}}\\[2ex]
&=&\int_1^{+\infty}\frac{-5/3+2H-1/3H^3}{\sqrt{2/3-H+H^3/3}\,H^{9/2}}\,dH.
\eeas
Let us make a further change of variables $t=1/H$ and integrate $I$ explicitly as follows.
\bes
I=\int_0^1\frac{-5/3t^4+2t^3-1/3}{\sqrt{2/3t^3-t^2+1/3}}\,dt=\int_0^1\frac{5t^3-t^2-t}{\sqrt{6t+3}}\,dt=\frac{1}{35(3 + \sqrt{3})}.
\lb{Ic}
\ees
\section{Connection between discrete and continuous \\coarsening laws} 
Here we show that the discrete coarsening law \rf{DL} can be recovered back from \rf{CRL} if the initial distribution $f(x)$ 
has the form \rf{DD} i.e. if it is represented by $k\in\N$ families as in \rf{hi} where we denote
\bes
\displaystyle i_m'=\lim_{N\go\infty}\frac{i_m}{N}.
\ees
In this case \rf{nh} implies
\bes
n(d)=1-\sum_{p=m}^ki_m'\ \ \text{if}\ \ d\in[d_m,d_{m-1}).
\ees
Substituting the last expression in \rf{CRL} implies
\beas
BT(d_m)&=&d_k\ln\left[\frac{1}{1-\sum_{p=m}^ki_m'}\right]+(d_{k-1}-d_k)\ln\left[\frac{1-i'_k}{1-\sum_{p=m}^ki_m'}\right]+...+\\
&+&(d_{m+1}-d_m)\ln\left[\frac{1-\sum_{p=m+1}^ki_m'}{1-\sum_{p=m}^ki_m'}\right]=BT(h_{m+1})+\\
&+&\left(Nd_m+\sum_{p=m}^{k}(d_p-d_m)i_p\right)\ln\left[\frac{1-\sum_{p=m+1}^ki_m'}{1-\sum_{p=m}^ki_m'}\right]
\eeas
Therefore,one obtains recursively that
\be
T(d_m)=\sum_{p=m}^{k}\frac{1}{B}\left(Nd_p+\sum_{p'=p}^{k}(d_p'-d_p)i_p'\right)\ln\left[\frac{1-\sum_{p=m+1}^ki_m'}{1-\sum_{p=m}^ki_m'}\right]
\lb{Tm3}
\ee
On the other hand dividing \rf{DL} by $N$ and proceeding to the limit $N\go\infty$  with $k$ fixed one obtains exactly \rf{Tm3}.
\end{appendix} 

\bibliographystyle{unsrtnat}
\bibliography{bibliography}
\clearpage
\addtocounter{tocdepth}{2}


\end{document}